\documentclass[a4paper]{article} 
\usepackage{amssymb}
\usepackage{amscd}
\usepackage{amsmath} 
\usepackage{graphicx}
\usepackage{latexsym} 
\usepackage{enumerate} 
\usepackage[usenames]{color}

\usepackage{theorem}
\theoremstyle{plain}
\theorembodyfont{\itshape}
\newtheorem{theorem}{Theorem}[section]
\newtheorem{proposition}[theorem]{Proposition}
\newtheorem{lemma}[theorem]{Lemma}
 
\theorembodyfont{\rmfamily}
\newtheorem{definition}[theorem]{Definition}
\newtheorem{remark}[theorem]{Remark}

\newtheorem{setup}[theorem]{Setup}

\makeatletter
\@addtoreset{figure}{section}
\@addtoreset{table}{section}
\makeatother
\def\bC{\mathbb{C}}
\def\bD{\mathbb{D}}

\def\bP{\mathbb{P}}
\def\bQ{\mathbb{Q}}
\def\bR{\mathbb{R}} 
\def\bZ{\mathbb{Z}}

\def\cC{\mathcal{C}}
\def\cE{\mathcal{E}} 

\def\cK{\mathcal{K}}
\def\cL{\mathcal{L}}

\def\cP{\mathcal{P}}
\def\cQ{\mathcal{Q}}
 
\def\rA{\mathrm{A}}

\def\Aut{\mathrm{Aut}} 
\def\rb{\mathrm{b}}

\def\rC{\mathrm{C}}
\def\ch{\mathrm{ch}}

\def\rD{\mathrm{D}}

\def\Dynkin{\mathrm{Dynkin}}
\def\rE{\mathrm{E}}

\def\Fix{\mathrm{Fix}}
\def\Fixe{\mathrm{Fix}^{\mathrm{e}}} 
\def\Fixi{\mathrm{Fix}^{\mathrm{i}}}

\def\GL{\mathrm{GL}}
\def\GR{\mathrm{GR}}

\def\rI{\mathrm{I}}

\def\rL{\mathrm{L}}

\def\ord{\mathrm{ord}}

\def\Pic{\mathrm{Pic}}

\def\Re{\mathrm{Re}}
\def\Res{\mathrm{Res}}
\def\rS{\mathrm{S}}
\def\Span{\mathrm{Span}} 
\def\ST{\mathrm{ST}}

\def\td{\mathrm{td}} 
\def\Tr{\mathrm{Tr}}
\def\rII{\mathrm{I\!I}}

\def\bA{\mbox{\boldmath $A$}}

\def\bB{\mbox{\boldmath $B$}}

\def\br{\mbox{\boldmath $r$}}

\def\bs{\mbox{\boldmath $s$}}

\def\bu{\mbox{\boldmath $u$}}

\def\bv{\mbox{\boldmath $v$}}

\def\b0{\mbox{\boldmath $0$}}

\def\d{\delta}
\def\vD{\varDelta}
\def\ve{\varepsilon}

\def\bAon{\bA_{\mathrm{on}}} 
\def\bAoff{\bA_{\mathrm{off}}}

\def\bBon{\bB_{\mathrm{on}}} 
 
\def\bBoff{\bB_{\mathrm{off}}}

\title{\bf K3 Surfaces, Picard Numbers and Siegel Disks}  
\author{Katsunori Iwasaki\thanks{Department of Mathematics, Faculty of Science, 
Hokkaido University, Kita 10, Nishi 8, Kita-ku, Sapporo 060-0810 Japan. 
{\tt iwasaki@math.sci.hokudai.ac.jp} (corresponding author).} \ and 
Yuta Takada\thanks{Department of Mathematics, Graduate School of Science, 
Hokkaido University, Kita 10, Nishi 8, Kita-ku, Sapporo 060-0810 Japan; 
JSPS Research Fellow. {\tt takada@math.sci.hokudai.ac.jp}}}
\date{November 6, 2021}   
\begin{document}
\maketitle
\begin{abstract} 
If a K3 surface admits an automorphism with a Siegel disk, then its Picard number 
is an even integer between $0$ and $18$.  
Conversely, using the method of hypergeometric groups, we are able to construct 
K3 surface automorphisms with Siegel disks that realize all possible Picard numbers. 
The constructions involve extensive computer searches for appropriate  
Salem numbers and computations of algebraic numbers arising from holomorphic  
Lefschetz-type formulas and related Grothendieck residues.     
\end{abstract} 
\section{Introduction} \label{sec:intro}
Let $X$ be a complex K3 surface, that is, a simply connected compact complex 
surface with trivial canonical bundle $K_X$.   
The middle cohomology group $H^2(X, \bZ)$ equipped with the intersection form 
is an even unimodular lattice of signature $(3, 19)$. 
The Hodge decomposition gives an orthogonal direct sum decomposition
$$
H^2(X, \bC) = H^{2,0}(X) \oplus H^{1,1}(X) \oplus H^{0,2}(X)
$$
of signatures $(1, 0) \oplus (1, 19) \oplus (1, 0)$. 
The Picard group (or N\'{e}ron-Severi group) of $X$ is the lattice 
$\Pic(X) = H^2(X, \bZ) \cap H^{1,1}(X)$, whose rank $\rho(X)$    
is called the {\sl Picard number} of $X$.  
It is an integer between $0$ and $20$.   
\par
Given a K3 surface automorphism $f : X \to X$, let $\lambda(f)$ be the spectral radius 
of $f^*|H^{1,1}(X)$. 
Then $\lambda(f) \ge 1$ and the topological entropy of $f$ is 
given by $h(f) = \log \lambda(f)$. 
There exists a constant $\delta(f) \in S^1$ such that 
$f^* \eta = \delta(f) \cdot \eta$ for a nowhere vanishing holomorphic 
$2$-form $\eta$ on $X$, where $S^1$ is the unit circle in $\bC$.   
Note that $f^*|H^{2, 0}(X)$ is the multiplication by $\delta(f)$. 
If $p \in X$ is a fixed point of $f$ then the holomorphic tangent map 
$(d f)_p : T_p X \to T_p X$ has determinant $\delta(f)$, 
so the number $\delta(f)$ is called the {\sl determinant} of $f$ by 
McMullen \cite{McMullen1}. 
It is referred to as the {\sl special eigenvalue} of $f$ in our previous 
paper \cite{IT}, where $\tau(f) := \delta(f) + \delta(f)^{-1}$ is called 
the {\sl special trace}.    
We remark that $\delta(f)$ is either a root of unity or a conjugate 
to a Salem number, and if $X$ is projective then $\delta(f)$ must be 
a root of unity. 
Here a {\sl Salem number} is an algebraic integer $\lambda > 1$ which is 
conjugate to $\lambda^{-1}$ and whose remaining conjugates lie on 
$S^1$.        
\par
Let $\bD$ be the unit disk in $\bC$.   
A map $R : (\bD^2, 0) \to (\bD^2, 0)$, $(z_1, z_2) \mapsto (\alpha_1 z_1, \alpha_2 z_2)$ 
with $\alpha_1$, $\alpha_2 \in S^1$ is said to be an {\sl irrational rotation} if 
$\alpha_1$ and $\alpha_2$ are multiplicatively independent, that is, if   
$\alpha_1^{m_1} \alpha_2^{m_2} = 1$ with $m_1$, $m_2 \in \bZ$ implies 
$m_1 = m_2 = 0$. 
Let $f : X \to X$ be an automorphism of a complex surface $X$. 
An open subset $U$ of $X$ is said to be a {\sl Siegel disk} for $f$ centered 
at $p \in U$ if $f$ preserves $(U, p)$ and $f|_{U} : (U, p) \to (U, p)$ is biholomorphically 
conjugate to an irrational rotation $R : (\bD^2, 0) \to (\bD^2, 0)$. 
If $X$ is a K3 surface and $f$ admits a Siegel disk, then $\lambda(f)$ must be 
a Salem number and $\delta(f)$ must be conjugate to $\lambda(f)$, in particular $X$ 
must be non-projective and $f$ must have a positive topological entropy 
(see McMullen \cite{McMullen1}).  
\par
McMullen \cite{McMullen1} synthesized examples of K3 surface automorphisms 
with a Siegel disk whose underlying K3 surfaces had Picard number $0$.   
Oguiso \cite{Oguiso} found an example of Picard number $8$. 
In \cite{IT} we constructed examples of Picard number $12$, whose entropy was 
the logarithm of Lehmer's number $\lambda_{\rL} \approx 1.17628$, the smallest 
Salem number ever known, as well as many more examples of Picard number $0$.  
The existence of a Siegel disk imposes a restriction on the Picard number of 
the underlying K3 surface.    
In this article we construct K3 surface automorphisms with Siegel disks that realize 
all possible Picard numbers. 
Our main result is stated as follows.     
\begin{theorem} \label{thm:main}
If a K3 surface $X$ admits an automorphism with at least one Siegel disk, 
then its Picard number $\rho(X)$ must be an even integer between $0$ and $18$. 
Conversely, for any such integer $r$ there exist K3 surface automorphisms 
$f : X \to X$ with Siegel disks such that $X$ has Picard number $\rho(X) = r$.  
\end{theorem} 
\par
The first half of the theorem is just a corollary to \cite[Theorem 7.4]{IT} and   
the essential part of the theorem is the second half stating that all Picard 
numbers $\rho = 0, 2, 4, \dots, 18$ can be realized by K3 surface automorphisms 
with Siegel disks.  
So this article is devoted to establishing the result in the second direction.    
\par
The construction of a K3 surface automorphism boils down 
to a lattice theoretic problem.     
Let $L$ be an abstract K3 lattice endowed with a Hodge structure  
$L_{\bC} = H^{2, 0} \oplus H^{1,1} \oplus H^{0,2}$ where $L_{\bC} := L \otimes \bC$.  
It determines the Picard lattice $\Pic := H^{1,1} \cap L$, root system  
$\vD := \{ \bu \in \Pic \mid (\bu, \bu) = -2 \}$ and Weyl group $W$, 
the group generated by reflections in root vectors. 
A positive cone $\cC^+$ is one of the two connected components 
of $\cC := \{ \bv \in H^{1,1}_{\bR} \mid (\bv, \bv) > 0 \}$, where 
$H^{1,1}_{\bR} :=  H^{1,1} \cap L_{\bR}$ with $L_{\bR} := L \otimes \bR$. 
We specify a Weyl chamber $\cK \subset \cC^+$ as the 
``K\"{a}hler cone''. 
This is equivalent to dividing $\vD$ into positive and negative 
roots $\vD = \vD^+ \amalg \vD^-$ in such a manner that   
$\cK = \{ \bv \in \cC^+ \mid \mbox{$(\bv, \bu) > 0$ for any  
$\bu \in \vD^+$} \}$. 
Note that $\vD^+$ determines a unique set of simple roots, 
say $\vD_{\rb}$, and vice versa. 
A Hodge isometry on $L$ is said to be positive if it preserves the 
connected components $\cC^{\pm}$ of $\cC$.  
It  falls into one of the three types; elliptic, parabolic and hyperbolic.  
By Torelli theorem and surjectivity of period mapping (see 
\cite[Chap.~V\!I\!I\!I]{BHPV}) any positive Hodge isometry 
$F : L \to L$ preserving the K\"{a}hler cone $\cK$ lifts to a unique K3 
surface automorphism $f : X \to X$ up to isomorphisms.  
In this article we deal with the case where $\Pic$ is negative definite, 
so that $\vD$ and $W$ are finite, and $F$ is a positive Hodge isometry 
of hyperbolic type. 
Then the resulting lift $f$ is a non-projective K3 surface automorphism 
of positive entropy.      
\par
We realize such structures by the method of hypergeometric 
groups developed in our article \cite{IT}. 
This method produces a large number of non-projective K3 surface 
automorphisms of positive entropy with various Picard numbers.    
From them we look for automorphisms with Siegel disks 
that cover all possible Picard numbers.  
\par
The plan of this article is as follows. 
In \S \ref{sec:hgm} we review our hypergeometric method in a way 
it is utilized in this article.  
In \S \ref{sec:al-cs}, by implementing the algorithm given in \S \ref{sec:hgm}, 
we develop extensive computer searches for pairs $(\varphi, \psi)$ leading to 
K3 surface automorphisms $f : X \to X$ of positive entropy with various 
Picard numbers $\rho$.  
The outputs are enormous, so only a part of which is exhibited in \S \ref{sec:al-cs}, 
with a more thorough presentation left to our web page \cite{IT2}. 
To pick out those entries with Siegel disks from so many candidates obtained 
in \S \ref{sec:al-cs}, we have to adjust Lefschetz-type fixed point 
formulas (FPF's) in a manner suitable for our purpose; to know how 
an isolated fixed point $p$ on the exceptional set $\cE(X)$ contributes  
to the FPF's; and to calculate the local index of $p$ as a Grothendieck residue, 
when $p$ is a multiple isolated fixed point.   
These tasks are done in \S \ref{sec:fpf}, \S \ref{sec:idx} and 
\S \ref{sec:res} respectively. 
In \S \ref{sec:SD}, combining all these ingredients with a criterion for 
Siegel disks, we construct K3 surface automorphisms with Siegel disks 
for Picard numbers $\rho = 2, 4, \dots, 18$ (see Theorem \ref{thm:SD}).  
The proofs for $\rho = 4, 6, \dots, 18$ are given in \S \ref{sec:SD}.  
The most difficult case of $\rho = 2$ is proved in a more general context 
in \S \ref{sec:P2} (see Theorem \ref{thm:P2}). 
The case of $\rho = 0$ is not treated in this article, as many examples are 
given in \cite{McMullen1} and \cite{IT}.               
\section{Method of Hypergeometric Groups} \label{sec:hgm}
To review the hypergeometric method, we recall some 
concepts and terminology on polynomials. 
In this section all polynomials are monic and defined over $\bZ$.    
Given a polynomial $u(z)$ of degree $n$, its {\sl reciprocal} is defined 
by $u^{\dagger}(z) := z^n \, u(z^{-1})$. 
We say that $u(z)$ is {\sl palindromic} if $u^{\dagger}(z) = u(z)$ and 
{\sl anti-palindromic} if $u^{\dagger}(z) = - u(z)$. 
If $u(z)$ is palindromic of even degree $n = 2m$, then there exists a unique 
polynomial $U(w)$ of degree $m$ such that $u(z) = z^m \, U(z+z^{-1})$.      
If $u(z)$ is anti-palindromic of even degree $n = 2m$, then there exists a unique 
polynomial $U(w)$ of degree $m-1$ such that $u(z) = (z-1)(z+1) z^{m-1} \, 
U(z+ z^{-1})$.  
In either case $U(w)$ is referred to as the {\sl trace polynomial} of $u(z)$.  
A palindromic polynomial $u(z)$ is said to be {\sl unramified} if 
$|u(\pm1)| = 1$. 
Such a polynomial is of degree even $n = 2 m$, has an even number, say $2 t$, 
of roots outside $S^1$ and satisfies  
\begin{equation} \label{eqn:sm}
t \equiv m \bmod 2, \qquad u(1) \cdot u(-1) = (-1)^m
\end{equation} 
(see Gross and McMullen \cite[Proposition 3.3]{GM}).  
The minimal polynomial of a Salem number is a {\sl Salem polynomial},  
which is palindromic of even degree and whose trace polynomial is called 
a {\sl Salem trace polynomial}.           
For any unramified Salem polynomial $u(z)$ the congruence in \eqref{eqn:sm}  
reads $m \equiv t = 1 \bmod 2$ and hence 
\begin{equation} \label{eqn:unrami}
\deg u(z) \equiv 2 \bmod 4. 
\end{equation} 
\par
Consider a coprime pair of anti-palindromic polynomial $\varphi(z)$ 
and palindromic polynomial $\psi(z)$ of degree $22$.  
Let $A$ and $B$ be the companion matrices of $\varphi(z)$ and $\psi(z)$ 
respectively, and let $H := \langle A, B \rangle \subset \GL(22, \bZ)$ 
be the hypergeometric group generated by $A$ and $B$. 
Then $C := A^{-1} B$ is a reflection, fixing a hyperplane in $\bQ^{22}$ 
pointwise and sending a nonzero vector $\br \in \bQ^{22}$ 
to its negative $-\br$.  
We have a free $\bZ$-module of rank $22$,  
\begin{equation} \label{eqn:L}
L = \langle \br, A \br, \dots, A^{21} \br \rangle_{\bZ} =  
\langle \br, B \br, \dots, B^{21} \br \rangle_{\bZ},  
\end{equation}
stable under the action of $H$. 
We can make $L$ into an $H$-invariant even lattice by providing it with the 
symmetric bilinear form $(A^{i-1} \br,  A^{j-1} \br) = \xi_{|i-j|}$, where 
$\xi_0 := 2$ and $\{ \xi_i \}_{i=1}^{\infty}$ is defined via the Taylor series expansion   
\begin{equation} \label{eqn:gram}
\dfrac{\psi(z)}{\varphi(z)} = 1 + \sum_{i=1}^{\infty} \xi_i \, z^{-i} 
\qquad \mbox{around} \quad z = \infty. 
\end{equation}
The Gram matrix $(B^{i-1} \br, B^{j-1} \br)$ for the $B$-basis is given by 
exchanging $\varphi(z)$ and $\psi(z)$ upside down in formula \eqref{eqn:gram}. 
The lattice $L$ is unimodular if and only if the resultant of 
$\varphi(z)$ and $\psi(z)$ satisfies 
\begin{equation} \label{eqn:um}
\Res(\varphi, \psi) = \pm 1, 
\end{equation}
in which case $\psi(z)$ must be unramified. 
Indeed, since $\varphi(z)$ is divisible by $(z-1)(z+1)$, the resultant is divisible 
by $\psi(1) \cdot \psi(-1)$ over $\bZ$, hence \eqref{eqn:um} implies 
$|\psi(\pm1)|=1$.  
For details we refer to \cite[Theorem 2.1]{IT}.    
\par
If the index of $L$ is positive, we replace $L$ by its negative $L(-1)$; 
otherwise, we keep $L$ as it is.    
This procedure is referred to as the {\sl renormalization} of $L$ and the 
renormalized bilinear form is called the {\sl intersection form} on $L$.  
In \cite{IT} we give a necessary and sufficient condition for the renormalized 
lattice $L$ to be a K3 lattice with a Hodge structure such that 
$A$ is a positive Hodge isometry of hyperbolic type. 
To review it, let $\Phi(w)$ and $\Psi(w)$ be the trace polynomials of 
$\varphi(z)$ and $\psi(z)$ respectively, that is, 
$$
\varphi(z) = (z-1)(z+1) z^{10} \, \Phi(z+ z^{-1}), \qquad 
\psi(z) = z^{11} \, \Psi(z+ z^{-1}). 
$$
Let $\bA$ be the multi-set of all complex roots of $\Phi(w)$ counted with 
multiplicity.  
Let $\bAon$ and $\bAoff$ be those parts of $\bA$ which lie on and off 
the interval $[-2, \, 2]$ respectively. 
Define $\bB$, $\bBon$ and $\bBoff$ in a similar manner for $\Psi(w)$. 
Then $\bAon$ and $\bBon$ dissect each other into interlacing components, 
called {\sl trace clusters}, such that 
\begin{equation} \label{eqn:trc}
-2 \le \bA_{s+1} < \bB_s < \bA_s < \dots < \bB_1 < \bA_1 \le 2, 
\end{equation}
where one or both of the end clusters $\bA_1$ and $\bA_{s+1}$ may be null, 
while any other cluster must be non-null.   
Put $\bA_{>2} := \bA \cap (2, \, \infty)$; $|\bAon|$ stands for the cardinality of 
$\bAon$ counted with multiplicity;  
$[\bAon] = 0^{\nu_0} 1^{\nu_1} 2^{\nu_2} 3^{\nu_3}$ means that 
$\bAon$ consists of $\nu_0$ null clusters, $\nu_1$ simple clusters, 
$\nu_2$ double clusters, $\nu_3$ triple clusters, where $j^{\nu_j}$ is  
omitted if $\nu_j = 0$.   
The same rule applies to $\bBon$ and other related entities.   
By ``doubles adjacent'' we mean the situation in which $\bAon$ and $\bBon$ 
contain unique double clusters $\bA_i$ and $\bB_j$ respectively, with $\bA_i$ 
and $\bB_j$ being adjacent to each other. 
If $\bA_i \cup \bB_j$ consists of four elements 
$x_1 < x_2 < x_3 < x_4$, then $x_2$ and $x_3$ are called the {\sl inner elements} 
of the adjacent pair (AP).   
As a part of \cite[Theorem 1.2]{IT} we have the following. 
\begin{table}[h]
\centerline{
\begin{tabular}{ccccccl|l}
\hline
\\[-4mm]
case & $s$ & $[\bAon]$ & $[\bBon]$ & $|\bA_{>2}|$ & $|\bBoff|$ & constraints & ST $\tau(A)$ 
\\[1mm]
\hline
\\[-4mm]
$1$ & $8$ & $0^21^63^1$ & $1^8$ & $1$ & $3$ & & middle of TC \\[1mm]
$2$ & $8$ & $0^21^63^1$ & $1^73^1$ & $1$ & $1$ & & middle of TC \\[1mm]
$3$ & $8$ & $0^11^72^1$ & $1^8$ & $1$ & $3$ & $|\bA_1| = 2$ & $\max \bA_1$ \\[1mm]
$4$ & $8$ & $0^11^72^1$ & $1^8$ & $1$ & $3$ & $|\bA_9|=2$ & $\min \bA_9$ \\[1mm]    
$5$ & $8$ & $0^11^72^1$ & $1^73^1$ & $1$ & $1$ & $|\bA_1|=2$ & $\max \bA_1$ \\[1mm] 
$6$ & $8$ & $0^11^72^1$ & $1^73^1$ & $1$ & $1$ & $|\bA_9|=2$ & $\min \bA_9$ \\[1mm] 
$7$ & $9$ & $0^21^72^1$ & $1^82^1$ & $1$ & $1$ & doubles adjacent & inner of AP \\[1mm] 
$8$ & $9$ & $0^11^9$ & $1^82^1$ & $1$ & $1$ & $|\bA_1|=1$, $|\bB_1|=2$ & element of $\bA_1$ \\[1mm] 
$9$ & $9$ & $0^11^9$ & $1^82^1$ & $1$ & $1$ & $|\bA_{10}|=1$, $|\bB_9|=2$ & element of $\bA_{10}$ \\[1mm]
\hline   
\end{tabular}} 
\caption{Conditions for $A$ to be a positive Hodge isometry of hyperbolic type \cite[Table 1.2]{IT}.} 
\label{tab:hyp-A} 
\end{table}
\begin{theorem} \label{thm:h-isom}
Let $L = L(\varphi, \psi)$ be a unimodular hypergeometric lattice of rank $22$. 
After renormalization, $L$ is a K3 lattice with a Hodge structure such that 
$A$ is a positive Hodge isometry of hyperbolic type, if and only if 
$\Phi(\pm2) \neq 0$, the roots of $\Phi(w)$ and $\Psi(w)$ are all simple 
and have any one of the configurations in Table $\ref{tab:hyp-A}$. 
In this case the special trace $\tau(A)$ and the Hodge structure up to complex 
conjugation are uniquely determined by the pair $(\varphi, \psi)$.   
The location of $\tau(A)$ is shown in the last column of Table $\ref{tab:hyp-A}$, 
where we mean by ``middle of TC" that $\tau(A)$ is 
the middle element of the unique triple cluster (TC) in $\bAon$, and 
by ``inner of AP" that $\tau(A)$ is the inner element in $\bAon$ of the unique 
AP of double clusters in $\bAon \cup \bBon$.  
\end{theorem}
\par
In the situation of Theorem \ref{thm:h-isom} $\varphi(z)$ factors as 
$\varphi(z) = \varphi_0(z) \cdot \varphi_1(z)$ where $\varphi_0(z)$ is a 
Salem polynomial and $\varphi_1(z)$ is a product of cyclotomic polynomials. 
Note that $\varphi_1(z)$ is divisible by $(z-1)(z+1)$. 
So we write 
\begin{equation} \label{eqn:SC}
\varphi_0(z) = S(z), \qquad \varphi_1(z) = (z-1)(z+1) \cdot C(z). 
\end{equation}
Let $\lambda(A) > 1$ be the Salem number associated with $S(z)$ and let 
$\delta(A)^{\pm1} \in S^1$ be the special eigenvalues corresponding to 
the special trace $\tau(A)$ in Theorem \ref{thm:h-isom}, that is, 
$\delta(A) + \delta(A)^{-1} = \tau(A)$. 
Then $\delta(A)$ is conjugate to $\lambda(A)$ and the Hodge structure (up to 
complex conjugation) is given by 
\begin{equation} \label{eqn:ell}
L_{\bC} = H^{2,0} \oplus H^{1,1} \oplus H^{0,2} 
= \ell \oplus (\ell \oplus \bar{\ell})^{\perp} \oplus \bar{\ell}, 
\end{equation}  
where $\ell$ is the eigen-line of $A$ corresponding to the eigenvalue $\delta(A)$ 
and $\bar{\ell}$ is the complex conjugate to $\ell$. 
Specify a positive cone $\cC^+ \subset H^{1,1}_{\bR}$ and put 
$\bs := S(A) \,\br$ with $\br$ being the vector in \eqref{eqn:L}.  
Then the intersection form is negative definite on the Picard lattice 
$\Pic := H^{1,1} \cap L$, whose rank, i.e. its Picard number is given by 
\begin{equation} \label{eqn:rho}
\rho = 22 - \deg S(z), 
\end{equation} 
and the vectors $\bs, A \bs, \dots, A^{\rho-1} \bs$ form a free basis, 
the {\sl standard basis}, of $\Pic$ (see \cite[Theorem 1.5]{IT}).  
The root system $\vD := \{ \bu \in \Pic : (\bu, \bu) = -2\}$ and the Weyl group 
$W$ are defined in the usual manner.  
The lexicographic order on $\Pic$ with respect to the standard basis leads to 
a set of positive roots $\vD^+$ and the corresponding Weyl chamber 
$\cK := \{ \bv \in \cC^+ : \mbox{$(\bv, \bu) > 0$ for any $\bu \in \vD^+$} \}$, 
which we specify as the ``K\"{a}hler cone".  
\par
The matrix $A$ may not preserve $\cK$, but there is a unique element 
$w_A \in W$ such that $\tilde{A} := w_A \circ A$ preserves $\cK$.  
We have an algorithm to determine $\vD$, $\vD^+$, $\vD_{\rb}$ and $w_A$ 
explicitly from the initial data $(\varphi, \psi)$, where $\vD_{\rb}$ is the 
simple system relative to $\vD^+$ (see \cite[Algorithm 7.5]{IT}). 
The Dynkin type of $\vD$ can be read off from the intersection relations for 
the simple roots in $\vD_{\rb}$.  
The characteristic polynomial $\tilde{\varphi}(z)$ of $\tilde{A}$ factors as 
\begin{equation} \label{eqn:char-tilde}
\tilde{\varphi}(z) = \varphi_0(z) \cdot \tilde{\varphi}_1(z), 
\end{equation}
where $\varphi_0(z) = S(z)$ is the same Salem polynomial as the one in \eqref{eqn:SC}  
while $\tilde{\varphi}_1(z)$ is a product of cyclotomic polynomial which, however, 
may differ from $\varphi_1(z)$ in \eqref{eqn:SC}. 
In particular $\tilde{A}$ and $A$ have the same spectral radius $\lambda(A)$ and 
the same special eigenvalue $\delta(A)$.  
Preserving the Hodge structure \eqref{eqn:ell} and the K\"{a}hler cone $\cK$, 
the modified matrix $\tilde{A}$ lifts to a K3 surface automorphism $f : X \to X$ 
of entropy $h(f) = \log \lambda(A)$ with special eigenvalue 
$\delta(f) = \delta(A)$, Picard lattice $\Pic(X) \cong \Pic$ and Picard number 
$\rho(X) = \rho$ given in \eqref{eqn:rho}. 
Recall that 
\begin{equation} \label{eqn:pic-char}
\mbox{$\tilde{\varphi}_1(z)$ is the characteristic polynomial of $f^*|\Pic(X)$}.       
\end{equation}    
Moreover, $\vD$, $\vD^+$ and $\vD_{\rb}$ lift to $\vD(X)$, $\vD^+(X)$ and 
$\vD_{\rb}(X)$ respectively, where $\vD(X)$ is the set of all $(-2)$-classes in 
$\Pic(X)$ with $\vD^+(X)$ being its subset of all effective $(-2)$-classes and 
$\vD_{\rb}(X)$ is the set of all $(-2)$-curves in $X$. 
How $f$ permutes the elements of $\vD_{\rb}(X)$ is faithfully represented 
by the action of $\tilde{A}$ on $\vD_{\rb}$.   
\section{Computer Searches} \label{sec:al-cs}
Let $\cP$ be a finite set of polynomials 
$\varphi(z) = (z-1)(z+1) \cdot S(z) \cdot C(z)$ of degree $22$
such that $S(z)$ is a Salem polynomial and $C(z)$ is 
a product of cyclotomic polynomials; see \eqref{eqn:SC}.  
Similarly let $\cQ$ be a finite set of unramified palindromic polynomials  
$\psi(z) \in \bZ[z]$ of degree $22$, where unramifiedness comes from 
the remark after \eqref{eqn:um}.   
For various choices of $\cP$ and $\cQ$ we make extensive computer searches 
for those pairs $(\varphi, \psi) \in \cP \times \cQ$ which satisfy firstly the 
unimodularity condition \eqref{eqn:um} and secondly all the conditions 
in Theorem \ref{thm:h-isom}.  
\par
The Salem numbers with any given degree, below any given bound, are finite in 
their cardinality. 
Thus we can speak of the $i$-th smallest Salem number $\lambda^{(d)}_i$ 
of degree $d$ and its minimal polynomial $\rS_i^{(d)}(z)$.  
The trace polynomial of $\rS_i^{(d)}(z)$ is denoted by $\ST_i^{(d)}(w)$.  
In his web page \cite{Mossinghoff} Mossinghoff gives a complete list of Salem 
numbers of small degrees, below certain bounds.      
A careful inspection of his tables together with the constraint 
\eqref{eqn:unrami} for unramifiedness leads us to the following observation.   
\begin{lemma} \label{lem:ursalem} 
Let $d$ be an even integer such that $4 \le d \le 22$. 
Then there exist exactly $N_d$ Salem numbers $\lambda$ of degree $d$ up to 
bound $\lambda \le M_d$, where $M_d$ and $N_d$ are given in 
Table $\ref{tab:ursalem}$ with $\mathrm{GR} := (1+ \sqrt{5})/2 \approx 1.61803$ 
being the golden ratio. 
Unramified Salem numbers of degree $d$ exist only when  
$d = 6$, $10$, $14$, $18$, $22$, for each of which there are exactly 
$N_d'$ such numbers up to bound $\lambda \le M_d$, where $N_d'$ is 
again given in Table $\ref{tab:ursalem}$.      
\begin{table}[h]
\centerline{
\begin{tabular}{c|cccccccccc|c}
\hline
           &       &           &          &          &           &           &           &            &           &         &     \\[-3mm]
$d$      & $4$ & $6$     & $8$   & $10$   & $12$   & $14$   & $16$   & $18$    & $20$   & $22$ &  total    \\[1mm]
\hline 
           &       &           &          &          &           &           &           &            &           &         &     \\[-3mm]
$M_d$  & $3$ & $2.8$  & $2.6$ & $2.4$  & $2.2$  & $2$     & $1.8$  & $\GR$ & $\GR$ & $1.5$ &  ---     \\[1mm]
$N_d$   &  $8$ & $34$  & $104$ & $223$ & $314$ & $390$ & $231$ & $141$  & $191$ &  $89$ & $1725$ \\[1mm]
$N_d'$  & --- & $3$    &  ---    & $29$  & ---     & $67$   &  ---   & $42$    & ---     & $30$  &  $171$ \\[1mm]
\hline
\end{tabular}}
\caption{Salem numbers $\lambda$ of degree $d \le 22$ up to bound $\lambda \le M_d$.}
\label{tab:ursalem}
\end{table}
\end{lemma}
\par
The $j$-th cyclotomic polynomial is denoted by $\rC_j(z)$. 
In \cite[\S5.2]{IT} we employ unconventional definitions $\rC_1(z) = (z-1)^2$ 
and $\rC_2(z) = (z+1)^2$ for $j = 1, 2$, but in this article we take the usual ones   
$\rC_1(z) = z-1$ and $\rC_2 = z+1$. 
For any $j \ge 3$ the polynomial $\rC_j(z)$ is palindromic of even degree, hence  
the congruence in \eqref{eqn:sm} with $t = 0$ implies that any unramified cyclotomic 
polynomial has a degree divisible by $4$.  
By \cite[Lemma 5.3]{IT} all unramified cyclotomic polynomials $\rC_l(z)$ of 
$\deg \rC_l(z) \le 16$ are exactly those with $l \in L_0$, where  
\begin{equation} \label{eqn:L0}
L_0 := \{12, 15, 20, 21, 24, 28, 30, 36, 40, 42, 48, 60\}.  
\end{equation}
\par
Seting up the ``principal" set $\cP$ is simple. 
Put $S(z) = \rS_i^{(d)}(z)$ for an even integer $d$ with $4 \le d \le 20$ and let  
\begin{equation} \label{eqn:C}
C(z) = \prod_{j \in J} \rC_j(z), \qquad \mbox{subject to the degree constraint} \quad 
d + \sum_{j \in J} \deg \rC_j(z) = 20,  
\end{equation}
where $J$ is a finite subset of $\bZ_{\ge 3}$ with $J = \emptyset$ for $d = 20$.    
We remark that $J$ contains neither $1$ nor $2$ because of $\Phi(\pm 2) \neq 0$ 
in Theorem \ref{thm:h-isom}. 
From \eqref{eqn:rho} the K3 surface to be constructed will have Picard number 
\begin{equation} \label{eqn:rho2}
\rho = 22-d.   \tag{$\ref{eqn:rho}'$}
\end{equation}  
Thus fixing a Picard number $\rho$ is fixing the degree $d$ according 
to \eqref{eqn:rho2} and specifying $\cP$ amounts to taking a finite subset of 
positive integers over  which the index $i$ of $\rS_i^{(d)}(z)$ ranges.    
The set of $J$'s is determined by $d$ according to \eqref{eqn:C}. 
In the Appendix we give a list of all Salem polynomials $\rS_i^{(d)}(w)$ that 
appear as $S(z)$ explicitly in this article; they are given in terms of their 
trace polynomials $\ST_i^{(d)}(w)$.     
\par
To set up the ``auxiliary" set $\cQ$ we observe from Table \ref{tab:hyp-A} that
the trace polynomial $\Psi(w)$ of $\psi(z)$ must have either ten or eight roots 
in $(-2, \, 2)$. 
An instance of the ten-root case is realized by the following setup.  
\begin{setup} \label{setup1} 
Let $\cQ$ be the set of all polynomials $\psi(z) = s(z) \cdot c(z)$ of degree $22$  
such that $s(z) = \rS^{(e)}_k(z)$ is an unramified Salem polynomial in 
Lemma \ref{lem:ursalem} and $c(z)$ is a product of unramified cyclotomic 
polynomials,  
\begin{equation} \label{eqn:c}
c(z) = \prod_{l \in L} \rC_l(z), \qquad \mbox{subject to the degree constraint} \quad 
e + \sum_{l \in L} \deg \rC_l(z) = 22.   
\end{equation}
Here since $e \ge 6$ it follows from \eqref{eqn:c} that $\deg \rC_l(z) \le 16$ 
for any $l \in L$, hence $L$ must be a subset of $L_0$ in \eqref{eqn:L0}.     
\end{setup}     
\par
For Picard numbers $\rho = 2, 4, 6, \dots, 16$, that is, for $d = 20, 18, 16, \dots, 6$, 
Setup \ref{setup1} with some choices of small indices $i$ (or even $i = 1$ only) for 
$S(z) = \rS_i^{(d)}(z)$ gives an abundance of solutions $(\varphi, \psi)$ satisfying 
the conditions in Theorem \ref{thm:h-isom}. 
We illustrate this by two computer outputs; one is for $\rho = 6$ $(d = 16)$, 
$i = 1, \dots,5$, and the other is for $\rho = 14$ $(d = 8)$, $i = 1, \dots, 16$. 
In these cases the results are given in Tables \ref{tab:P6} and \ref{tab:P14} respectively. 
In Table \ref{tab:P14} almost all solutions with $i = 2, \dots, 15$ are omitted because 
there are too many of them. 
We refer to our web page \cite{IT2} for more extensive outputs that cover all of  
the cases $\rho = 2, 4, 6, \dots, 16$. 
\begin{table}[p]
\centerline{
\begin{tabular}{ccccccccc}
\hline
          &           &           &           &       &          &                                &                        &      \\[-3mm]
$S(z)$ & $C(z)$ & $s(z)$ & $c(z)$ & ST & Dynkin & $\Tilde{\varphi}_1(z)$ & $\Tr \Tilde{A}$ & SD \\ 
\hline
$\rS^{(16)}_{1}$ & $\rC_8$ & $\rS^{(10)}_{4}$ & $\rC_{36}$ & $\tau_2$ & $\rE_6$ & $\rC_1^4\rC_2^2$ & $3$ \\
$\rS^{(16)}_{1}$ & $\rC_8$ & $\rS^{(10)}_{5}$ & $\rC_{21}$ & $\tau_2$ & $\rE_6$ & $\rC_1^4\rC_2^2$ & $3$ \\
$\rS^{(16)}_{1}$ & $\rC_8$ & $\rS^{(10)}_{7}$ & $\rC_{28}$ & $\tau_5$ & $\rE_6$ & $\rC_1^4\rC_2^2$ & $3$ \\
$\rS^{(16)}_{1}$ & $\rC_8$ & $\rS^{(14)}_{120}$ & $\rC_{15}$ & $\tau_1$ & $\rE_6$ & $\rC_1^4\rC_2^2$ & $3$ \\
$\rS^{(16)}_{1}$ & $\rC_8$ & $\rS^{(22)}_{7}$ & $1$ & $\tau_6$ & $\rE_6$ & $\rC_1^4\rC_2^2$ & $3$ \\
$\rS^{(16)}_{1}$ & $\rC_8$ & $\rS^{(22)}_{85}$ & $1$ & $\tau_5$ & $\rE_6$ & $\rC_1^4\rC_2^2$ & $3$ \\
$\rS^{(16)}_{1}$ & $\rC_{10}$ & $\rS^{(10)}_{5}$ & $\rC_{21}$ & $\tau_5$ & $\rE_6$ & $\rC_1^4\rC_2^2$ & $3$ \\
$\rS^{(16)}_{1}$ & $\rC_{10}$ & $\rS^{(10)}_{15}$ & $\rC_{28}$ & $\tau_1$ & $\rE_6$ & $\rC_1^4\rC_2^2$ & $3$ \\
$\rS^{(16)}_{1}$ & $\rC_{10}$ & $\rS^{(14)}_{4}$ & $\rC_{15}$ & $\tau_2$ & $\rE_6$ & $\rC_1^4\rC_2^2$ & $3$ \\
$\rS^{(16)}_{1}$ & $\rC_{10}$ & $\rS^{(14)}_{4}$ & $\rC_{24}$ & $\tau_2$ & $\rE_6$ & $\rC_1^4\rC_2^2$ & $3$ \\
$\rS^{(16)}_{1}$ & $\rC_{10}$ & $\rS^{(14)}_{15}$ & $\rC_{15}$ & $\tau_1$ & $\rE_6$ & $\rC_1^4\rC_2^2$ & $3$ \\
$\rS^{(16)}_{1}$ & $\rC_{10}$ & $\rS^{(14)}_{15}$ & $\rC_{24}$ & $\tau_1$ & $\rE_6$ & $\rC_1^4\rC_2^2$ & $3$ \\
$\rS^{(16)}_{1}$ & $\rC_{10}$ & $\rS^{(22)}_{85}$ & $1$ & $\tau_2$ & $\rE_6$ & $\rC_1^4\rC_2^2$ & $3$ \\
$\rS^{(16)}_{1}$ & $\rC_3\rC_4$ & $\rS^{(22)}_{41}$ & $1$ & $\tau_5$ & $\rE_6$ & $\rC_1^4\rC_2^2$ & $3$ \\ 
\hline
$\rS^{(16)}_{2}$ & $\rC_5$ & $\rS^{(22)}_{18}$ & $1$ & $\tau_5$ & $\rA_6$ & $\rC_1^6$ & $5$ \\
$\rS^{(16)}_{2}$ & $\rC_3\rC_4$ & $\rS^{(22)}_{18}$ & $1$ & $\tau_3$ & $\rA_6$ & $\rC_1^6$ & $5$ \\
$\rS^{(16)}_{2}$ & $\rC_3\rC_6$ & $\rS^{(22)}_{1}$ & $1$ & $\tau_1$ & $\rA_6$ & $\rC_1^6$ & $5$ \\
$\rS^{(16)}_{2}$ & $\rC_3\rC_6$ & $\rS^{(22)}_{18}$ & $1$ & $\tau_2$ & $\rA_6$ & $\rC_1^6$ & $5$ \\
$\rS^{(16)}_{2}$ & $\rC_3\rC_6$ & $\rS^{(22)}_{72}$ & $1$ & $\tau_4$ & $\rA_6$ & $\rC_1^6$ & $5$ \\ 
\hline
$\rS^{(16)}_{3}$ & $\rC_{10}$ & $\rS^{(18)}_{1}$ & $\rC_{12}$ & $\tau_1$ & $\rA_4$ & $\rC_1^3\rC_2^3$ & $0$ \\
$\rS^{(16)}_{3}$ & $\rC_{10}$ & $\rS^{(22)}_{2}$ & $1$ & $\tau_1$ & $\rA_4$ & $\rC_1^3\rC_2^3$ & $0$ \\
$\rS^{(16)}_{3}$ & $\rC_3\rC_4$ & $\rS^{(22)}_{6}$ & $1$ & $\tau_3$ & $\rA_2$ & $\rC_1^3\rC_2\rC_4$ & $2$ \\ 
\hline
$\rS^{(16)}_{4}$ & $\rC_5$ & $\rS^{(22)}_{22}$ & $1$ & $\tau_6$ & $\rD_6$ & $\rC_1^5\rC_2$ & $5$ \\
$\rS^{(16)}_{4}$ & $\rC_{10}$ & $\rS^{(10)}_{15}$ & $\rC_{28}$ & $\tau_4$ & $\rD_6$ & $\rC_1^5\rC_2$ & $5$ \\
$\rS^{(16)}_{4}$ & $\rC_{10}$ & $\rS^{(14)}_{4}$ & $\rC_{15}$ & $\tau_4$ & $\rD_6$ & $\rC_1^5\rC_2$ & $5$ \\
$\rS^{(16)}_{4}$ & $\rC_{10}$ & $\rS^{(22)}_{13}$ & $1$ & $\tau_3$ & $\rD_6$ & $\rC_1^5\rC_2$ & $5$ \\
$\rS^{(16)}_{4}$ & $\rC_{10}$ & $\rS^{(22)}_{32}$ & $1$ & $\tau_1$ & $\rD_6$ & $\rC_1^5\rC_2$ & $5$ \\
$\rS^{(16)}_{4}$ & $\rC_{10}$ & $\rS^{(22)}_{39}$ & $1$ & $\tau_5$ & $\rD_6$ & $\rC_1^5\rC_2$ & $5$ \\
$\rS^{(16)}_{4}$ & $\rC_{10}$ & $\rS^{(22)}_{85}$ & $1$ & $\tau_5$ & $\rD_6$ & $\rC_1^5\rC_2$ & $5$ \\ 
\hline
$\rS^{(16)}_{5}$ & $\rC_5$ & $\rS^{(14)}_{3}$ & $\rC_{30}$ & $\tau_3$ & $\rD_6$ & $\rC_1^5\rC_2$ & $4$ \\
$\rS^{(16)}_{5}$ & $\rC_5$ & $\rS^{(14)}_{12}$ & $\rC_{30}$ & $\tau_6$ & $\rD_6$ & $\rC_1^5\rC_2$ & $4$ \\
$\rS^{(16)}_{5}$ & $\rC_5$ & $\rS^{(22)}_{10}$ & $1$ & $\tau_5$ & $\rD_6$ & $\rC_1^5\rC_2$ & $4$ \\
$\rS^{(16)}_{5}$ & $\rC_{10}$ & $\rS^{(10)}_{12}$ & $\rC_{36}$ & $\tau_4$ & $\rD_6$ & $\rC_1^5\rC_2$ & $4$ \\
$\rS^{(16)}_{5}$ & $\rC_{10}$ & $\rS^{(22)}_{2}$ & $1$ & $\tau_4$ & $\rD_6$ & $\rC_1^5\rC_2$ & $4$ \\
$\rS^{(16)}_{5}$ & $\rC_{10}$ & $\rS^{(22)}_{3}$ & $1$ & $\tau_4$ & $\rD_6$ & $\rC_1^5\rC_2$ & $4$ \\
$\rS^{(16)}_{5}$ & $\rC_{10}$ & $\rS^{(22)}_{10}$ & $1$ & $\tau_2$ & $\rD_6$ & $\rC_1^5\rC_2$ & $4$ \\
$\rS^{(16)}_{5}$ & $\rC_{12}$ & $\rS^{(10)}_{24}$ & $\rC_{42}$ & $\tau_7$ & $\rA_1^{\oplus 2}\oplus \rD_4$ & $\rC_1^3\rC_2\rC_3$ & $1$ \\
$\rS^{(16)}_{5}$ & $\rC_{12}$ & $\rS^{(10)}_{107}$ & $\rC_{42}$ & $\tau_1$ & $\rA_1^{\oplus 2}\oplus \rD_4$ & $\rC_1^3\rC_2\rC_3$ & $1$ \\
$\rS^{(16)}_{5}$ & $\rC_{12}$ & $\rS^{(14)}_{12}$ & $\rC_{30}$ & $\tau_2$ & $\rA_1^{\oplus 2}\oplus \rD_4$ & $\rC_1^3\rC_2\rC_3$ & $1$ \\
$\rS^{(16)}_{5}$ & $\rC_{12}$ & $\rS^{(14)}_{17}$ & $\rC_{15}$ & $\tau_4$ & $\rA_1^{\oplus 2}\oplus \rD_4$ & $\rC_1^3\rC_2\rC_3$ & $1$ \\
$\rS^{(16)}_{5}$ & $\rC_{12}$ & $\rS^{(14)}_{382}$ & $\rC_{20}$ & $\tau_1$ & $\rA_1^{\oplus 2}\oplus \rD_4$ & $\rC_1^3\rC_2\rC_3$ & $1$ \\
$\rS^{(16)}_{5}$ & $\rC_3\rC_4$ & $\rS^{(10)}_{2}$ & $\rC_{42}$ & $\tau_6$ & $\emptyset$ & $\rC_1\rC_2\rC_3\rC_4$ & $-1$ & S \\
$\rS^{(16)}_{5}$ & $\rC_3\rC_4$ & $\rS^{(14)}_{11}$ & $\rC_{30}$ & $\tau_1$ & $\emptyset$ & $\rC_1\rC_2\rC_3\rC_4$ & $-1$ & S \\
$\rS^{(16)}_{5}$ & $\rC_3\rC_4$ & $\rS^{(22)}_{10}$ & $1$ & $\tau_3$ & $\emptyset$ & $\rC_1\rC_2\rC_3\rC_4$ & $-1$ & S \\
$\rS^{(16)}_{5}$ & $\rC_3\rC_4$ & $\rS^{(22)}_{43}$ & $1$ & $\tau_4$ & $\emptyset$ & $\rC_1\rC_2\rC_3\rC_4$ & $-1$ & S \\ 
\hline
\end{tabular}} 
\caption{Picard number $\rho = 6$ (Setup \ref{setup1}).} 
\label{tab:P6}
\end{table}
\begin{table}[p]
\centerline{
\begin{tabular}{ccccccccc}
\hline
          &           &           &          &      &           &                               &                       &      \\[-3mm]
$S(z)$ & $C(z)$ & $s(z)$ & $c(z)$ & ST & Dynkin & $\tilde{\varphi}_1(t)$ & $\Tr \tilde{A}$ & SD \\ 
\hline
$\rS^{(8)}_{1}$ & $\rC_{28}$ & $\rS^{(10)}_{24}$  & $\rC_{36}$ & $\tau_2$ & $\rA_2$ & $\rC_1 \rC_2 \rC_{28}$ & $0$ \\
$\rS^{(8)}_{1}$ & $\rC_{28}$ & $\rS^{(14)}_{3}$    & $\rC_{30}$ & $\tau_1$ & $\rA_2$ & $\rC_1 \rC_2 \rC_{28}$ & $0$ \\
$\rS^{(8)}_{1}$ & $\rC_{28}$ & $\rS^{(14)}_{234}$ & $\rC_{30}$ & $\tau_3$ & $\rA_2$ & $\rC_1 \rC_2 \rC_{28}$ & $0$ \\
$\rS^{(8)}_{1}$ & $\rC_{28}$ & $\rS^{(18)}_{65}$ & $\rC_{12}$ & $\tau_2$ & $\rA_2$ & $\rC_1 \rC_2 \rC_{28}$ & $0$ \\
$\rS^{(8)}_{1}$ & $\rC_{28}$ & $\rS^{(18)}_{109}$ & $\rC_{12}$ & $\tau_3$ & $\rA_2$ & $\rC_1 \rC_2 \rC_{28}$ & $0$ \\
$\rS^{(8)}_{1}$ & $\rC_4\rC_{11}$ & $\rS^{(6)}_{1}$ & $\rC_{60}$ & $\tau_1$ & $\rA_1^{\oplus 2}$ & $\rC_1^2\rC_2^2\rC_{11}$ & $-1$ \\
$\rS^{(8)}_{1}$ & $\rC_4\rC_{11}$ & $\rS^{(10)}_{2}$ & $\rC_{42}$ & $\tau_2$ & $\rA_1^{\oplus 2}$ & $\rC_1^2\rC_2^2\rC_{11}$ & $-1$ \\
$\rS^{(8)}_{1}$ & $\rC_5\rC_{20}$ & $\rS^{(10)}_{1}$ & $\rC_{42}$ & $\tau_1$ & $\rA_1^{\oplus 5}$ & $\rC_1\rC_2\rC_5\rC_{20}$ & $-1$ \\
$\rS^{(8)}_{1}$ & $\rC_5\rC_{20}$ & $\rS^{(10)}_{1}$ & $\rC_{12}\rC_{24}$ & $\tau_2$ & $\rA_1^{\oplus 5}$ & $\rC_1\rC_2\rC_5\rC_{20}$ & $-1$ \\
$\rS^{(8)}_{1}$ & $\rC_5\rC_{20}$ & $\rS^{(10)}_{2}$ & $\rC_{12}\rC_{30}$ & $\tau_1$ & $\rA_1^{\oplus 5}$ & $\rC_1\rC_2\rC_5\rC_{20}$ & $-1$ \\
$\rS^{(8)}_{1}$ & $\rC_5\rC_{20}$ & $\rS^{(14)}_{7}$ & $\rC_{30}$ & $\tau_1$ & $\rA_1^{\oplus 5}$ & $\rC_1\rC_2\rC_5\rC_{20}$ & $-1$ \\
$\rS^{(8)}_{1}$ & $\rC_5\rC_{20}$ & $\rS^{(18)}_{48}$ & $\rC_{12}$ & $\tau_2$ & $\rA_1^{\oplus 5}$ & $\rC_1\rC_2\rC_5\rC_{20}$ & $-1$ \\
$\rS^{(8)}_{1}$ & $\rC_5\rC_{20}$ & $\rS^{(18)}_{95}$ & $\rC_{12}$ & $\tau_1$ & $\rA_1^{\oplus 5}$ & $\rC_1\rC_2\rC_5\rC_{20}$ & $-1$ \\
$\rS^{(8)}_{1}$ & $\rC_5\rC_{20}$ & $\rS^{(22)}_{18}$ & $1$ & $\tau_1$ & $\rA_1^{\oplus 5}$ & $\rC_1\rC_2\rC_5\rC_{20}$ & $-1$ \\
$\rS^{(8)}_{1}$ & $\rC_{10}\rC_{15}$ & $\rS^{(10)}_{7}$ & $\rC_{12}\rC_{24}$ & $\tau_2$ & $\rE_6$ & $\rC_1^4\rC_2^2\rC_{15}$ & $3$ \\
$\rS^{(8)}_{1}$ & $\rC_{10}\rC_{15}$ & $\rS^{(10)}_{15}$ & $\rC_{42}$ & $\tau_2$ & $\rE_6$ & $\rC_1^4\rC_2^2\rC_{15}$ & $3$ \\
$\rS^{(8)}_{1}$ & $\rC_{10}\rC_{15}$ & $\rS^{(18)}_{4}$ & $\rC_{12}$ & $\tau_2$ & $\rE_6$ & $\rC_1^4\rC_2^2\rC_{15}$ & $3$ \\
$\rS^{(8)}_{1}$ & $\rC_{10}\rC_{15}$ & $\rS^{(18)}_{7}$ & $\rC_{12}$ & $\tau_2$ & $\rE_6$ & $\rC_1^4\rC_2^2\rC_{15}$ & $3$ \\
$\rS^{(8)}_{1}$ & $\rC_{10}\rC_{15}$ & $\rS^{(18)}_{33}$ & $\rC_{12}$ & $\tau_2$ & $\rE_6$ & $\rC_1^4\rC_2^2\rC_{15}$ & $3$ \\
$\rS^{(8)}_{1}$ & $\rC_{10}\rC_{15}$ & $\rS^{(18)}_{43}$ & $\rC_{12}$ & $\tau_3$ & $\rE_6$ & $\rC_1^4\rC_2^2\rC_{15}$ & $3$ \\
$\rS^{(8)}_{1}$ & $\rC_{10}\rC_{15}$ & $\rS^{(22)}_{3}$ & $1$ & $\tau_2$ & $\rE_6$ & $\rC_1^4\rC_2^2\rC_{15}$ & $3$ \\
$\rS^{(8)}_{1}$ & $\rC_{10}\rC_{15}$ & $\rS^{(22)}_{13}$ & $1$ & $\tau_2$ & $\rE_6$ & $\rC_1^4\rC_2^2\rC_{15}$ & $3$ \\
$\rS^{(8)}_{1}$ & $\rC_{10}\rC_{15}$ & $\rS^{(22)}_{32}$ & $1$ & $\tau_2$ & $\rE_6$ & $\rC_1^4\rC_2^2\rC_{15}$ & $3$ \\
$\rS^{(8)}_{1}$ & $\rC_{10}\rC_{24}$ & $\rS^{(22)}_{10}$ & $1$ & $\tau_2$ & $\rE_6\oplus \rE_8$ & $\rC_1^{12}\rC_2^2$ & $10$ \\
$\rS^{(8)}_{1}$ & $\rC_3\rC_4\rC_{30}$ & $\rS^{(22)}_{72}$ & $1$ & $\tau_1$ & $\rE_6\oplus \rE_8$ & $\rC_1^{12}\rC_2^2$ & $10$ \\
$\rS^{(8)}_{1}$ & $\rC_3\rC_{12}\rC_{18}$ & $\rS^{(10)}_{12}$ & $\rC_{42}$ & $\tau_3$ & $\rA_2^{\oplus 3}\oplus \rE_7$ & $\rC_1^8\rC_2^2\rC_3\rC_6$ & $6$ \\
$\rS^{(8)}_{1}$ & $\rC_3\rC_{12}\rC_{18}$ & $\rS^{(14)}_{3}$ & $\rC_{30}$ & $\tau_3$ & $\rA_2^{\oplus 3}\oplus \rE_7$ & $\rC_1^8\rC_2^2\rC_3\rC_6$ & $6$ \\
$\rS^{(8)}_{1}$ & $\rC_3\rC_{12}\rC_{18}$ & $\rS^{(22)}_{52}$ & $1$ & $\tau_3$ & $\rA_2^{\oplus 3}\oplus \rE_7$ & $\rC_1^8\rC_2^2\rC_3\rC_6$ & $6$ \\
$\rS^{(8)}_{1}$ & $\rC_4\rC_7\rC_8$ & $\rS^{(22)}_{72}$ & $1$ & $\tau_2$ & $\rA_7$ & $\rC_1^7\rC_2\rC_4\rC_8$ & $6$ \\ 
\hline
$\rS_2^{(8)}$ & $\rC_{13}$ & $\rS_1^{(10)}$ & $\rC_{21}$ & $\tau_3$ & $\emptyset$ & $\rC_1 \rC_2 \rC_{13}$ & $0$ \\
$\rS_2^{(8)}$ & $\rC_{13}$ & $\rS_1^{(10)}$ & $\rC_{36}$ & $\tau_2$ & $\emptyset$ & $\rC_1 \rC_2 \rC_{13}$ & $0$ \\
$\vdots$ & $\vdots$ & $\vdots$   & $\vdots$  & $\vdots$  & $\vdots$   &  $\vdots$  & $\vdots$  &   \\
$\rS_{15}^{(8)}$ & $\rC_7 \rC_9$ & $\rS_2^{(14)}$ & $\rC_{15}$ & $\tau_3$ & $\rE_6$ & $\rC_1^7 \rC_2 \rC_7$ & $7$ \\  
$\rS_{15}^{(8)}$ & $\rC_8 \rC_{15}$ & $\rS_{145}^{(10)}$ & $\rC_{42}$ & $\tau_1$ & $\rE_8$ & $\rC_1^9 \rC_2 \rC_8$ & $10$ \\ 
$\rS_{15}^{(8)}$ & $\rC_8 \rC_{15}$ & $\rS_{145}^{(10)}$ & $\rC_{12} \rC_{20}$ & $\tau_2$ & $\rE_8$ & $\rC_1^9 \rC_2 \rC_8$ & $10$ \\ 
\hline
$\rS^{(8)}_{16}$ & $\rC_{13}$ & $\rS^{(10)}_{2}$ & $\rC_{42}$ & $\tau_3$ & $\mathrm{A}_1\oplus \mathrm{A}_{13}$ & $\rC_1^{14}$ & $14$ \\
$\rS^{(8)}_{16}$ & $\rC_3\rC_4\rC_{20}$ & $\rS^{(10)}_{2}$ & $\rC_{42}$ & $\tau_2$ & $\rA_1^{\oplus 10}\oplus \rA_2$ & $\rC_1^4\rC_2^2\rC_5\rC_{10}$ & $2$ \\
$\rS^{(8)}_{16}$ & $\rC_3\rC_{12}\rC_{18}$ & $\rS^{(10)}_{12}$ & $\rC_{42}$ & $\tau_2$ & $\emptyset$ & $\rC_1 \rC_2 \rC_3 \rC_{12} \rC_{18}$ & $-1$ & S \\ 
\hline
\end{tabular}}
\caption{Picard number $\rho = 14$ (Setup \ref{setup1}).} 
\label{tab:P14}
\par\vspace{10mm}\noindent
\centerline{
\begin{tabular}{ccccccccc}
\hline
           &          &           &           &      &           &                             &                      &      \\[-3mm]
$S(z)$  & $C(z)$ & $s(z)$  & $c(z)$ & ST & Dynkin & $\tilde{\varphi}_1(z)$ & $\Tr \tilde{A}$ & SD \\ 
\hline
$\rS^{(4)}_{1}$ & $\rC_3\rC_{18}\rC_{24}$ & $\rS^{(14)}_{48}$ & $\rC_{30}$ & $\tau_1$ & $\rA_1^{\oplus 3}\oplus \rA_2^{\oplus 3}$ & $\rC_1^2\rC_2^2 \rC_3^2 \rC_6 \rC_{24}$ & $0$ \\ 
\hline
\end{tabular}}
\caption{Picard number $\rho = 18$ (Setup \ref{setup1}).}
\label{tab:P18}
\end{table}
\begin{table}[p]
\centerline{
\begin{tabular}{cccccccc}
\hline
          &           &               &      &                 &                               &                       &      \\[-3mm]          
$S(z)$ & $C(z)$ & $\psi(z)$ & ST & $\Dynkin$ & $\tilde{\varphi}_1(z)$ & $\Tr \tilde{A}$ & SD \\
\hline
$\rS_1^{(4)}$ & $\rC_{17}$ & $579$ & $\tau_1$ & $\emptyset$ & $\rC_{17}$ & $0$ &  \\
$\rS_1^{(4)}$ & $\rC_{32}$ & $289$ & $\tau_1$ & $\emptyset$ & $\rC_{32}$ & $1$ &  \\
$\rS_1^{(4)}$ & $\rC_{32}$ & $576$ & $\tau_1$ & $\emptyset$ & $\rC_{32}$ & $1$ &  \\
$\rS_1^{(4)}$ & $\rC_{32}$ & $692$ & $\tau_1$ & $\emptyset$ & $\rC_{32}$ & $1$ &  \\
$\rS_1^{(4)}$ & $ \rC_{32}$ & $711$ & $\tau_1$ & $\emptyset$ & $\rC_{32}$ & $1$ &  \\
$\rS_1^{(4)}$ & $ \rC_{40}$ & $40$ & $\tau_1$ & $\rA_2$ & $\rC_1 \rC_2 \rC_{40}$ & $1$ &  \\
$\rS_1^{(4)}$ & $ \rC_{40}$ & $58$ & $\tau_1$ & $\rA_2$ & $\rC_1 \rC_2 \rC_{40}$ & $1$ &  \\
$\rS_1^{(4)}$ & $ \rC_{40}$ & $515$ & $\tau_1$ & $\rA_2$ & $\rC_1 \rC_2 \rC_{40}$ & $1$ &  \\
$\rS_1^{(4)}$ & $ \rC_{40}$ & $579$ & $\tau_1$ & $\rA_2$ & $\rC_1 \rC_2 \rC_{40}$ & $1$ &  \\
$\rS_1^{(4)}$ & $ \rC_{40}$ & $873$ & $\tau_1$ & $\rA_2$ & $\rC_1 \rC_2 \rC_{40}$ & $1$ &  \\
$\rS_1^{(4)}$ & $ \rC_{48}$ & $692$ & $\tau_1$ & $\rA_2$ & $\rC_1 \rC_2 \rC_{48}$ & $1$ &  \\
$\rS_1^{(4)}$ & $ \rC_{48}$ & $699$ & $\tau_1$ & $\rA_2$ & $\rC_1 \rC_2 \rC_{48}$ & $1$ &  \\
$\rS_1^{(4)}$ & $ \rC_{60}$ & $457$ & $\tau_1$ & $\rA_2$ & $\rC_1 \rC_2 \rC_{60}$ & $1$ &  \\
$\rS_1^{(4)}$ & $ \rC_{60}$ & $699$ & $\tau_1$ & $\rA_2$ & $\rC_1 \rC_2 \rC_{60}$ & $1$ &  \\
$\rS_1^{(4)}$ & $ \rC_{60}$ & $744$ & $\tau_1$ & $\rA_2$ & $\rC_1 \rC_2 \rC_{60}$ & $1$ &  \\
$\rS_1^{(4)}$ & $ \rC_{60}$ & $961$ & $\tau_1$ & $\rA_2$ & $\rC_1 \rC_2 \rC_{60}$ & $1$ &  \\
$\rS_1^{(4)}$ & $ \rC_5 \rC_{26}$ & $664$ & $\tau_1$ & $\rA_1^{\oplus 5}$ & $\rC_1 \rC_2 \rC_5 \rC_{26}$ & $1$ &  \\
$\rS_1^{(4)}$ & $ \rC_5 \rC_{26}$ & $679$ & $\tau_1$ & $\rA_1^{\oplus 5}$ & $\rC_1 \rC_2 \rC_5 \rC_{26}$ & $1$ &  \\
$\rS_1^{(4)}$ & $ \rC_5 \rC_{26}$ & $792$ & $\tau_1$ & $\rA_1^{\oplus 5}$ & $\rC_1 \rC_2 \rC_5 \rC_{26}$ & $1$ &  \\
$\rS_1^{(4)}$ & $ \rC_5 \rC_{26}$ & $893$ & $\tau_1$ & $\rA_1^{\oplus 5}$ & $\rC_1 \rC_2 \rC_5 \rC_{26}$ & $1$ &  \\
$\rS_1^{(4)}$ & $ \rC_5 \rC_{26}$ & $961$ & $\tau_1$ & $\rA_1^{\oplus 5}$ & $\rC_1 \rC_2 \rC_5 \rC_{26}$ & $1$ &  \\
$\rS_1^{(4)}$ & $ \rC_5 \rC_{36}$ & $873$ & $\tau_1$ & $\rA_1^{\oplus 5}\oplus \rE_6^{\oplus 2}$ & $\rC_1^5 \rC_2^5 \rC_4^2 \rC_5$ & $0$ &  \\
$\rS_1^{(4)}$ & $ \rC_5 \rC_{36}$ & $901$ & $\tau_1$ & $\rA_1^{\oplus 5}\oplus \rE_6^{\oplus 2}$ & $\rC_1^5 \rC_2^5 \rC_4^2 \rC_5$ & $0$ &  \\
$\rS_1^{(4)}$ & $ \rC_5 \rC_{36}$ & $961$ & $\tau_1$ & $\rA_1^{\oplus 5}\oplus \rE_6^{\oplus 2}$ & $\rC_1^5 \rC_2^5 \rC_4^2 \rC_5$ & $0$ &  \\
$\rS_1^{(4)}$ & $ \rC_8 \rC_{36}$ & $457$ & $\tau_1$ & $\rE_6^{\oplus 3}$ & $\rC_1^8 \rC_2^6 \rC_4^2$ & $3$ &  \\
$\rS_1^{(4)}$ & $ \rC_8 \rC_{36}$ & $515$ & $\tau_1$ & $\rE_6^{\oplus 3}$ & $\rC_1^8 \rC_2^6 \rC_4^2$ & $3$ &  \\
$\rS_1^{(4)}$ & $ \rC_8 \rC_{36}$ & $699$ & $\tau_1$ & $\rE_6^{\oplus 3}$ & $\rC_1^8 \rC_2^6 \rC_4^2$ & $3$ &  \\
$\rS_1^{(4)}$ & $ \rC_8 \rC_{36}$ & $712$ & $\tau_1$ & $\rE_6^{\oplus 3}$ & $\rC_1^8 \rC_2^6 \rC_4^2$ & $3$ &  \\
$\rS_1^{(4)}$ & $ \rC_3 \rC_4 \rC_{28}$ & $515$ & $\tau_1$ & $\rA_3$ & $\rC_1^3 \rC_2 \rC_4 \rC_{28}$ & $3$ &  \\
$\rS_1^{(4)}$ & $ \rC_3 \rC_4 \rC_{28}$ & $870$ & $\tau_1$ & $\rA_3$ & $\rC_1^3 \rC_2 \rC_4 \rC_{28}$ & $3$ &  \\
$\rS_1^{(4)}$ & $\rC_3 \rC_9 \rC_{15}$ & $870$ & $\tau_1$ & $\emptyset$ & $\rC_3 \rC_9 \rC_{15}$ & $1$ &  \\
$\rS_1^{(4)}$ & $ \rC_3 \rC_9 \rC_{24}$ & $692$ & $\tau_1$ & $\emptyset$ & $\rC_3 \rC_9 \rC_{24}$ & $0$ &  \\
$\rS_1^{(4)}$ & $ \rC_4 \rC_{10} \rC_{11}$ & $692$ & $\tau_1$ & $\rD_{11}$ & $\rC_1^{11} \rC_2 \rC_4 \rC_{10}$ & $12$ &  \\
$\rS_1^{(4)}$ & $ \rC_8 \rC_{12} \rC_{30}$ & $523$ & $\tau_1$ & $\rA_2^{\oplus 2}\oplus \rE_6\oplus \rE_8$ & $\rC_1^{13} \rC_2^3 \rC_4$ & $11$ & S \\
$\rS_1^{(4)}$ & $ \rC_9 \rC_{10} \rC_{18}$ & $711$ & $\tau_1$ & $\emptyset$ & $\rC_9 \rC_{10} \rC_{18}$ & $2$ &  \\
$\rS_1^{(4)}$ & $ \rC_{10} \rC_{12} \rC_{16}$ & $259$ & $\tau_1$ & $\rA_2 \oplus \rA_2$ & $\rC_1^2 \rC_2^2 \rC_4 \rC_{10} \rC_{16}$ & $2$ &  \\
$\rS_1^{(4)}$ & $ \rC_3 \rC_4 \rC_6 \rC_{11}$ & $279$ & $\tau_1$ & $\rA_{11}$ & $\rC_1^{11} \rC_2 \rC_3 \rC_4 \rC_6$ & $11$ &  \\
$\rS_1^{(4)}$ & $ \rC_3 \rC_4 \rC_8 \rC_{15}$ & $515$ & $\tau_1$ & $\rA_2^{\oplus 4}$ & $\rC_1^2 \rC_2^2 \rC_4^2 \rC_3 \rC_{15}$ & $1$ &  \\
$\rS_1^{(4)}$ & $ \rC_3 \rC_4 \rC_8 \rC_{15}$ & $870$ & $\tau_1$ & $\rA_2^{\oplus 4}$ & $\rC_1^2 \rC_2^2 \rC_4^2 \rC_3 \rC_{15}$ & $1$ &  \\
$\rS_1^{(4)}$ & $\rC_3 \rC_4 \rC_8 \rC_{16}$ & $699$ & $\tau_1$ & $\rA_2$ & $\rC_1^3 \rC_2 \rC_4 \rC_8 \rC_{16}$ & $3$ &  \\
\hline
\end{tabular}}	
\caption{Picard number $\rho = 18$ (Setup \ref{setup2}).} 
\label{tab:setup2}
\end{table}
\begin{remark} \label{rem:howto} 
We explain how to look at Tables \ref{tab:P6}, \ref{tab:P14} and  
similar tables to be given later.    
The meanings of the $S(z)$, $C(z)$, $s(z)$, $c(z)$ columns are clear.  
The ST column indicates the value of the special trace $\tau$, where 
\begin{equation} \label{eqn:trace}
\tau_0 >  \tau_1  > \tau_2 > \cdots > \tau_{d/2-1} \qquad \mbox{with} \quad 
\tau_0 > 2 > \tau_1     
\end{equation}
are the roots of the trace polynomial $ST(w)$ of $S(z)$.  
Notice that $\tau_j$ depends on the index $(d, i)$ of $ST(w) = \ST_i^{(d)}(w)$, 
so should be written $\tau_{i,j}^{(d)}$ to be precise, but this dependence is 
suppressed in notation.     
The Dynkin column exhibits the Dynkin type of the root system $\vD$. 
The $\tilde{\varphi}_1(z)$ column shows the $\tilde{\varphi}_1(z)$-component 
of the characteristic polynomial $\tilde{\varphi}(z)$ of  
$\tilde{A}$; see \eqref{eqn:char-tilde}. 
In \S \ref{sec:SD} Lefschetz-type fixed point formulas will be used to 
look for Siegel disks and the values of $\Tr \, \tilde{A}$ are needed to do so, 
hence this information is given in the $\Tr \, \tilde{A}$ column.  
If a solution is marked with S in the SD column, then it can be shown that the  
automorphism arising from this entry has at least one Siegel disk. 
(But a blank in this column does not claim non-existence of Siegel disks.)  
\end{remark}
\par
For $\rho = 18$, i.e. $d = 4$, however, Setup \ref{setup1} with 
$S(z) = \rS_1^{(4)}(z)$ leads to only one solution in Table \ref{tab:P18}, for which 
it is difficult to decide whether the resulting automorphism has Siegel disks or not. 
Thus we propose an alternative setup which offers a wider variety of candidates for 
$\psi(z)$.    
\begin{setup} \label{setup2}
For Picard number $\rho = 18$, take $S(z) := \rS^{(4)}_1(z) = 
z^4 - z^3 -z^2 - z + 1$ to be the minimal polynomial of the Salem 
number $\lambda^{(4)}_1 \approx 1.7220838$, 
and let $\cQ$ be the set of all unramified palindromic polynomials,   
$$
\psi(z) = z^{22} + c_1 z^{21} + \cdots + c_{10} z^{12} + c_{11} z^{11} + c_{10} z^{10} + \cdots 
+ c_1 z + 1 \in \bZ[z]  
$$ 
such that the following three conditions are satisfied:   
\begin{enumerate}
\setlength{\itemsep}{-1pt}
\item[(i)] $c_j \in \{0, \pm 1, \pm 2\}$ for $j = 1, \dots, 9$,   
\item[(ii)] the trace polynomial $\Psi(w)$ of $\psi(z)$ has ten or eight roots on the 
interval $(-2, \, 2)$, 
\item[(iii)] a part of unimodularity condition: the resultant of $\rS_1^{(4)}(z)$ and  
$\psi(z)$ is $\pm 1$. 
\end{enumerate} 
\end{setup}
\par
Unramifiedness of $\psi(z)$ implies $\psi(1) = \pm 1$ and $\psi(-1) = \mp 1$, 
where $\psi(\pm1)$ must have different signs because the second formula 
in \eqref{eqn:sm} gives $\psi(1) \cdot \psi(-1) = -1$. 
Thus $c_{10}$ and $c_{11}$ can be determined from $(c_1, \dots, c_9)$ by  
$$
c_{10} = -1-c_2-c_4-c_6-c_8, \qquad c_{11} = c_{11}^{\pm} := \pm 1-2(c_1+c_3+c_5+c_7+c_9),   
$$
in a unique manner for $c_{10}$ and in two ways for $c_{11}$.  
A computer enumeration shows that $\cQ$ contains a total of $1019$ polynomials.  
They can be identified by the numbering according to the lexicographical order for  
words $(c_1, \dots, c_{11})$.  
All solutions to Setup \ref{setup2} are given in Table \ref{tab:setup2}, where $\psi(z)$ is 
shown by its ID number and $\tau_1 := (1-\sqrt{13})/2 \approx -1.30278$
 is the only root in $(-2, \, 2)$ of $\ST^{(4)}_1(w) = w^2-w-3$.  
For the entry marked with S in the SD column, $\psi(z)$ has ID number $523$. 
Explicitly, this polynomial is given by 
\begin{equation} \label{eqn:psi18} 
\psi(z) 
= z^{22} -z^{21} -2 z^{20} +2 z^{18} + z^{17} -z^{15} -2 z^{14} + z^{12} +z^{11} 
+ z^{10} -2 z^8 - z^7 + z^5 +2 z^4 -2 z^2 -z +1.    
\end{equation}
It is the minimal polynomial of a Salem number $\lambda \approx 1.72654$ 
of degree $22$, which does not appear in Mossinghoff's list  
\cite{Mossinghoff} because $\lambda$ is beyond his bound $M_{22} = 1.5$ in 
Table \ref{tab:ursalem}. 
It is why Setup \ref{setup1} fails to find this solution.          
\section{Fixed Point Formulas} \label{sec:fpf}
We present two fixed point formulas (FPF's), originally due to Saito \cite{Saito}, 
Toledo and Tong \cite{TT}, which are needed to discuss the existence of Siegel disks. 
Let $f : X \to X$ be a K3 surface automorphism such that  
\begin{enumerate} 
\setlength{\itemsep}{-1pt}
\item[(C1)] $X$ is non-projective and the intersection form on 
$\Pic(X)$ is negative definite,   
\item[(C2)] the special eigenvalue $\delta = \delta(f)$ is conjugate 
to a Salem number.
\end{enumerate} 
These conditions are satisfied by all non-projective K3 surface automorphisms  
produced by the method of hypergeometric groups \cite[Theorem 1.5]{IT}. 
Looking for Siegel disks naturally involves questions about fixed points of $f$.  
We have to control the fixed point set of $f$, which consists of 
{\sl isolated} fixed points and possibly occurring {\sl fixed curves}. 
{\sl Invariant} (but not fixed) {\sl curves} should also be relevant to this issue.    
By condition (C1) any irreducible curve in $X$ is a $(-2)$-curve 
\cite[Lemma 7.3]{IT}. 
This fact and triviality of the canonical bundle $K_X$ are helpful in 
discussing questions about fixed curves and invariant curves.  
We begin by recalling the following.  
\begin{lemma} \label{lem:transv} 
If two distinct $(-2)$-curves in $X$ meet then they 
meet exactly in one point transversally. 
\end{lemma}       
{\it Proof}. If $C_1$ and $C_2$ are such curves, then $1 \le C_1 \cdot C_2$ and 
$(C_1+C_2)^2 = 2 C_1 \cdot C_2 + C_1^2 + C_2^2 = 2(C_1 \cdot C_2 -2) \le - 2$, 
since $C_1$ and $C_2$ are $(-2)$-curves, the intersection form on $\Pic(X)$ is 
even and negative definite, and $C_1 + C_2 \neq 0$ in $\Pic(X)$. 
Therefore $C_1 \cdot C_2 = 1$ and the assertion follows. 
\hfill $\Box$ \par\medskip
A fixed point $p \in X$ of $f$ is {\sl isolated} if and only if its {\sl multiplicity}
\begin{equation} \label{eqn:index}
\mu_p(f) := \dim_{\bC} \, \left( \bC\{ z \}/\frak{a} \right) \qquad 
\mbox{with} \quad \frak{a} := (z_1 - f_1(z), \, z_2-f_2(z)), 
\end{equation}
is finite,  where $(f_1, f_2)$ is the local representation of $f$ in terms of 
a local chart $z = (z_1, z_2)$ around $p \leftrightarrow z = (0, 0)$, 
$\bC\{ z \}$ is the convergent power series ring in two variables $z = (z_1, z_2)$  
and $\frak{a}$ is its ideal generated by $z_1-f_1(z)$ and $z_2-f_2(z)$.  
Let $\Fixi(f)$ denote the set of all isolated fixed point of $f$.  
\begin{proposition} \label{prop:saito} 
If $N_f$ is the number of $(-2)$-curves fixed pointwise by $f$, then    
\begin{equation} \label{eqn:saito}
\sum_{p \in \Fixi(f)} \mu_p(f) = \Tr \, f^*|H^2(X, \bC) + 2(1- N_f).   
\end{equation}
\end{proposition}
{\it Proof}. 
We use S.~Saito's fixed point formula \cite[formula (0.2)]{Saito} which is stated as    
$$
L(f) := \sum_{j=0}^4 (-1)^j \, \Tr \, f^*|H^j(X, \bC)   
= \sum_{p \in X_0(f)} \mu_p(f) + \sum_{C \in X_{\rI}(f)} \chi_{C} \cdot \mu_C(f) + 
\sum_{C \in X_{\rII}(f)} \tau_{C} \cdot \mu_C(f), 
$$
where $X_0(f)$ is the set of all fixed points of $f$ while $X_{\rI}(f)$ and $X_{\rII}(f)$ 
are the sets of all irreducible fixed curves of types I and I\!I respectively, 
$\chi_C$ is the Euler number of the normalization of $C$ and $\tau_C$ is the 
self-intersection number of $C$. 
For the definitions of Saito's indices $\mu_p(f)$ and $\mu_C(f)$ 
we refer to \cite[\S 3]{IU}. 
For any isolated fixed point $p$, Saito's index $\mu_p(f)$ coincides with the 
the multiplicity defined in \eqref{eqn:index}, hence the same notation is 
employed for the two concepts.    
The formula holds for compact K\"{a}hler surfaces \cite[Theorem 4.3]{DNT}. 
Since $X$ is a K3 surface we have $L(f) = 2 + \Tr \, f^*|H^2(X)$. 
Any fixed curve $C$ is a $(-2)$-curve isomorphic to $\bP^1$. 
The differential $d f$ acts on the normal bundle $N_C$ to $C$ 
as multiplication by $\delta \neq 1$.   
Thus $C \in X_{\rI}(f)$ with $\chi_C = 2$ and $X_{\rII}(f)$ is empty.  
An inspection shows that $\mu_C(f) = 1$ and $\mu_p(f) = 0$ at each  
$p \in C$ (see also \cite[\S9.3]{IT}). 
Putting all these facts into Saito's formula we obtain formula \eqref{eqn:saito}. 
\hfill $\Box$ \par\medskip 
The holomorphic local index of an isolated fixed point $p \in \Fixi(f)$ is 
given by the Grothendieck residue  
\begin{equation} \label{eqn:g-res}
\nu_p(f) = \Res_p \, \omega \qquad \mbox{with} \quad 
\omega := \dfrac{d z_1 \wedge d z_2}{(z_1 - f_1(z)) (z_2 - f_2(z))}.  
\end{equation}
If $p$ is simple i.e. $\mu_p(f) = 1$ or equivalently if $p$ is transverse to the effect 
that the tangent map $(d f)_p$ does not have eigenvalue $1$, then the index 
$\nu_p(f)$ admits a simpler representation 
\begin{equation} \label{eqn:h-index}
\nu_p(f) = \frac{1}{\det(I - (d f)_p)} = \dfrac{1}{1- \Tr (d f)_p + \delta}.   
\end{equation}
\begin{proposition} \label{prop:TT} 
If $\delta = \delta(f)$ is the special eigenvalue of $f$ and 
$N_f$ is the number in Proposition $\ref{prop:saito}$, then 
\begin{equation} \label{eqn:TT} 
1 + \delta^{-1} = \sum_{p \in \Fixi(f)} \nu_p(f) + N_f \dfrac{1+\delta}{(1 - \delta)^2}.  
\end{equation}
\end{proposition}
{\it Proof}. 
We use the Toledo-Tong fixed point formula \cite[Theorem (4.10)]{TT} in 
$2$-dimensional case. 
If any isolated fixed point $p \in X$ is transverse and if any connected 
component $C$ of the $1$-dimensional fixed point set is also transverse  
to the effect that $C$ is a smooth curve and the induced differential map 
$d^N \! f$ on the normal line bundle $N_C$ to $C$ has eigenvalue 
$\lambda_C \neq 1$, then the holomorphic Lefschetz number $\cL(f)$ 
is expressed as     
\begin{equation} \label{eqn:TT2}
\cL(f) := \sum_{j=0}^2 (-1)^j \, \Tr \, f^*|H^{0, j}(X)  
= \sum_{p \in \Fixi(f)} \nu_p(f) + \sum_{C} \nu_C(f).  
\end{equation}
Here $\nu_p(f)$ is given by \eqref{eqn:h-index}, while if $TC$ is the tangent bundle 
to $C$ and $\check{N}_C$ is the dual bundle to $N_C$, then 
$$
\nu_C(f) = \int_C \dfrac{ \td(TC)}{1 - \lambda_C \cdot \ch(\check{N}_C)}  
= \dfrac{1}{1- \lambda_C} \int_C 
\left\{ 
\dfrac{1}{2} \, c_1(TC) + \dfrac{\lambda_C \cdot c_1(\check{N}_C)}{1-\lambda_C} 
\right\},  
$$
where $\td$ and $\ch$ stand for Todd class and Chern character respectively. 
When $p \in \Fixi(f)$ is not transverse, $\nu_p(f)$ can be expressed  
by the Grothendieck residue in \eqref{eqn:g-res} 
(see Toledo \cite[formula (6.3)]{Toledo}).  
\par
Currently, we have $\cL(f) = 1+\bar{\delta} = 1+ \delta^{-1}$ since 
$H^{0,1}(X) = 0$, $H^{0,2}(X) = \overline{H^{2,0}(X)}$ and $f^*|H^{2, 0}(X) = \delta$. 
Let $C$ be any connected component of the $1$-dimensional 
fixed point set. 
If $C$ contains two distinct $(-2)$-curves  $C_1$ and $C_2$ meeting in 
a point $p$, then Lemma \ref{lem:transv} shows that they meets  
transversally in $p$, so $(d f)_p$ acts on $T_p X = T_p C_1 \oplus T_p C_2$ 
trivially, but this contradicts the fact that $\det (d f)_p = \delta \neq 1$.    
Therefore $C$ is just a single $(-2)$-curve, which is smooth. 
Triviality of the canonical bundle $K_X$ implies that $\check{N}_C$ is 
isomorphic to the tangent bundle $TC$.  
Taking $C \cong \bP^1$, $\lambda_C = \delta$ and 
$\int_C c_1(T \bP^1) = 2$ into account, we have  
\begin{equation} \label{eqn:mu_C}
\nu_C(f) = \dfrac{1+\delta}{(1-\delta)^2} 
\end{equation}
for any $(-2)$-curve $C$ fixed by $f$. 
Thus Toledo-Tong formula \eqref{eqn:TT2} leads to the equation 
\eqref{eqn:TT}. \hfill $\Box$  
\section{Indices on Exceptional Set} \label{sec:idx}
The union $\cE = \cE(X)$ of all $(-2)$-curves in $X$ is referred to as 
the {\sl exceptional set}.  
We are interested in how the isolated fixed points on $\cE$ contribute to the 
FPF's \eqref{eqn:saito} and \eqref{eqn:TT}.  
This problem may be considered component-wise for each connected component 
$\cE'$ of $\cE$ preserved by $f$. 
In what follows we denote by $\mu(f, \cE')$ and $\nu(f, \cE')$ the sum of 
$\mu_p(f)$ and that of $\nu_p(f)$ taken over all isolated fixed points $p$ on 
$\cE'$ respectively. 
\par
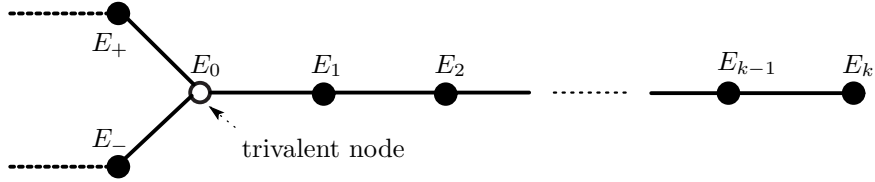
\begin{figure}[h] 
\centerline{ 
{\unitlength 0.1in%
\begin{picture}(44.1800,9.0200)(18.4000,-18.4100)%
%
\special{pn 20}%
\special{ar 2823 1403 51 51 0.0000000 6.2831853}%
%
\special{sh 1.000}%
\special{ia 3463 1411 51 51 0.0000000 6.2831853}%
\special{pn 20}%
\special{ar 3463 1411 51 51 0.0000000 6.2831853}%
%
\special{sh 1.000}%
\special{ia 4095 1411 51 51 0.0000000 6.2831853}%
\special{pn 20}%
\special{ar 4095 1411 51 51 0.0000000 6.2831853}%
%
\special{sh 1.000}%
\special{ia 5559 1406 51 51 0.0000000 6.2831853}%
\special{pn 20}%
\special{ar 5559 1406 51 51 0.0000000 6.2831853}%
%
\special{sh 1.000}%
\special{ia 6207 1406 51 51 0.0000000 6.2831853}%
\special{pn 20}%
\special{ar 6207 1406 51 51 0.0000000 6.2831853}%
%
\special{pn 20}%
\special{pa 2879 1403}%
\special{pa 4527 1403}%
\special{fp}%
%
\special{pn 20}%
\special{pa 5159 1406}%
\special{pa 6223 1406}%
\special{fp}%
\special{pa 6223 1406}%
\special{pa 6223 1406}%
\special{fp}%
%
\special{pn 13}%
\special{pa 4655 1403}%
\special{pa 5055 1403}%
\special{dt 0.045}%
\put(27.7000,-13.2000){\makebox(0,0)[lb]{$E_0$}}%
\put(33.9600,-13.2000){\makebox(0,0)[lb]{$E_1$}}%
\put(40.3600,-13.2000){\makebox(0,0)[lb]{$E_2$}}%
\put(55.0000,-13.2000){\makebox(0,0)[lb]{$E_{k-1}$}}%
\put(61.4800,-13.2000){\makebox(0,0)[lb]{$E_k$}}%
%
\special{pn 8}%
\special{pa 3004 1584}%
\special{pa 2876 1480}%
\special{dt 0.045}%
\special{sh 1}%
\special{pa 2876 1480}%
\special{pa 2915 1538}%
\special{pa 2917 1514}%
\special{pa 2940 1507}%
\special{pa 2876 1480}%
\special{fp}%
\special{pa 2876 1480}%
\special{pa 2876 1480}%
\special{dt 0.045}%
%
\special{sh 1.000}%
\special{ia 2400 990 51 51 0.0000000 6.2831853}%
\special{pn 20}%
\special{ar 2400 990 51 51 0.0000000 6.2831853}%
%
\special{sh 1.000}%
\special{ia 2400 1790 51 51 0.0000000 6.2831853}%
\special{pn 20}%
\special{ar 2400 1790 51 51 0.0000000 6.2831853}%
%
\special{pn 20}%
\special{pa 2420 990}%
\special{pa 2790 1360}%
\special{fp}%
%
\special{pn 20}%
\special{pa 2400 1780}%
\special{pa 2780 1420}%
\special{fp}%
%
\special{pn 20}%
\special{pa 1840 990}%
\special{pa 2380 990}%
\special{dt 0.045}%
%
\special{pn 20}%
\special{pa 1840 1790}%
\special{pa 2380 1790}%
\special{dt 0.045}%
\put(22.6000,-12.2000){\makebox(0,0)[lb]{$E_+$}}%
\put(22.6000,-17.2000){\makebox(0,0)[lb]{$E_-$}}%
\put(30.4000,-17.4000){\makebox(0,0)[lb]{trivalent node}}%
\end{picture}}%
 }
\caption{Dynkin diagram $\varGamma$ with a trivalent node $E_0$.}
\label{fig:Dynkin}
\end{figure}
We discuss the case where the dual graph of $\cE'$ is a Dynkin diagram 
$\varGamma$ with a {\sl trivalent} node, that is, of type $\rD$ or $\rE$ 
as in Figure \ref{fig:Dynkin}.  
If $\varGamma$ is of type $\rE_7$ or $\rE_8$ then the automorphism group 
$\Aut \, \varGamma$ is trivial, while if $\varGamma$ is of type $\rD_n$ $(n \ge 5)$ 
or $\rE_6$ then $\Aut \, \varGamma \cong \bZ/2\bZ$, where the nontrivial 
automorphism fixes an arm $E_1, \dots, E_k$ emanating from the trivalent node $E_0$, 
but permutes the remaining two arms, namely, those containing $E_{\pm}$.         
\begin{lemma} \label{lem:DE} 
Let $\cE'$ be a connected component of $\cE$ preserved by $f$, the  
dual graph of which is a Dynkin diagram $\varGamma$ with a trivalent node. 
Then all isolated fixed points $p \in \cE'$ are simple $\mu_p(f) = 1$, that is, 
transverse.      
\begin{enumerate}
\setlength{\itemsep}{-1pt}
\item If $f$ acts on $\varGamma$ trivially, then $\cE'$ contains exactly  
one irreducible fixed curve, $\mu(f, \cE') = n-1$ and   
\begin{subequations} \label{eqn:nup}
\begin{alignat}{2}
\nu(f, \cE') 
&= -\dfrac{\delta}{(1-\delta)^2} \left( \dfrac{2}{1+\delta} 
+ \dfrac{1+ \delta + \cdots + \delta^{n-4}}{1+\delta+ \cdots +\delta^{n-3}}\right)  
& \qquad & \mbox{for type $\rD_n$, $n \ge 5$}, \label{eqn:nupD} \\
\nu(f, \cE') &=-\dfrac{\delta}{(1-\delta)^2} \left( \dfrac{1}{1+\delta} 
+ \dfrac{1+\delta}{1+\delta+\delta^2} 
+ \dfrac{1+\delta + \cdots + \delta^{n-5}}{1+\delta + \cdots +\delta^{n-4}}\right)  
& \qquad & \mbox{for type $\rE_n$, $n = 6, 7, 8$}. \label{eqn:nupE}
\end{alignat}
\end{subequations}
\item If $f$ acts on $\varGamma$ non-trivially, then $\cE'$ contains no 
irreducible fixed curve and 
\begin{subequations} \label{eqn:num}
\begin{alignat}{3} 
\mu(f, \cE') &= n-1, \qquad &
\nu(f, \cE') &= \dfrac{1}{2(1+\delta)} + 
\dfrac{1 + \delta + \cdots + \delta^{n-3}}{2(1+\delta^{n-2})}  
& \qquad & \mbox{for type $\rD_n$, $n \ge 5$}, \label{eqn:numD} \\
\mu(f, \cE') &= 3, \qquad &
\nu(f, \cE') &= \dfrac{1}{2(1+\delta)} + 
\dfrac{1+\delta}{2(1+\delta^2)} & \qquad & \mbox{for type $\rE_6$}. 
\label{eqn:numE} 
\end{alignat}
\end{subequations}
\end{enumerate}
\end{lemma}
{\it Proof}.  
Let $q_j$ be the intersection of $E_j$ and $E_{j+1}$ for $j = 0, \dots, k-1$. 
Let $q_{\pm}$ be the intersection of $E_0$ and $E_{\pm}$. 
\par
\begin{figure}[h]
\centerline{\input{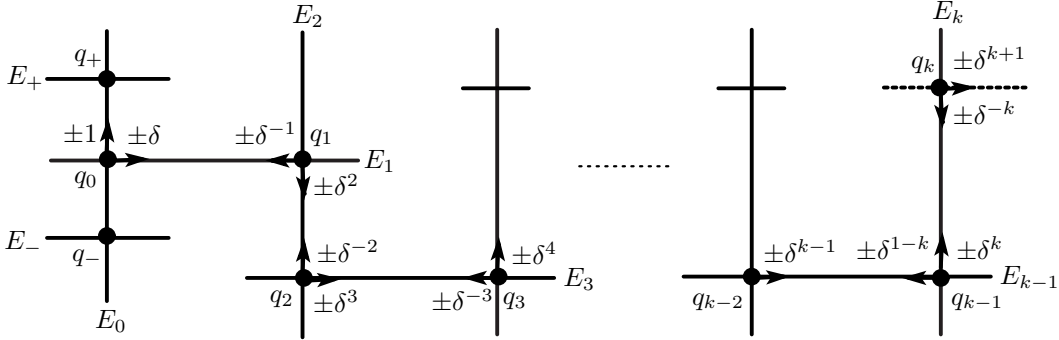}}
\caption{An arm $E_1 \cup \cdots \cup E_k$ emanating from the curve $E_0$ of a trivalent node.}
\label{fig:excep} 
\end{figure}

Assertion (1). 
In this case $E_0$ is a fixed curve of $f$, since the M\"{o}bius transformation $f_{E_0}$ 
fixes the three points $q_0$ and $q_{\pm}$. 
Thus one has $(d f_{E_0})_{q_0} = 1$, $(d f_{E_1})_{q_0} = \delta$, and continues to 
get $(d f_{E_j})_{q_j} = \delta^{-j}$ and $(d f_{E_{j+1}})_{q_j} = \delta^{j+1}$ 
successively for $j = 1, \dots, k-1$, until arriving at a unique fixed point 
$q_k \in E_k$ different from $q_{k-1}$, at which $(d f)_{q_k}$ has eigenvalues 
$\delta^{-k}$ and $\delta^{k+1}$ (see Figure \ref{fig:excep}).   
This argument works smoothly and shows that $q_1, \dots, q_k$ are transverse 
fixed points, because $\delta$ is not a root of unity.   
Moreover, those $k$ points are all of the isolated fixed points on 
$E_0 \cup \cdots \cup E_k$, and the sum of their $\nu$-indices can be calculated as    
$$
\varLambda_k^+ := \sum_{j=1}^k \nu_{q_j}(f) = 
\sum_{j=1}^k \dfrac{1}{(1-\delta^{-j})(1-\delta^{j+1})} 
= - \dfrac{\delta(1+\delta + \cdots + \delta^{k-1})}{(1-\delta)^2(1 + \delta + \cdots + \delta^k)}.
$$
If $\varGamma$ is of type $\rD_n$, then the three arms of $\varGamma$ have lengths 
$1$, $1$ and $n-3$, so $\mu(f, \cE') = 1+1+(n-3) = n-1$ and 
$\nu(f, \cE') = \varLambda_1^+ +\varLambda_1^+ + \varLambda_{n-3}^+$; 
this yields \eqref{eqn:nupD}. 
If $\varGamma$ is of type $\rE_n$ $(n = 6, 7, 8)$, then the three arms of $\varGamma$ 
have lengths $1$, $2$ and $n-4$, so $\mu(f, \cE') = 1+2+(n-4) = n-1$ and 
$\nu(f, \cE') = \varLambda_1^+ + \varLambda_2^+ + \varLambda_{n-3}^+$; 
this yields \eqref{eqn:nupE}. 
\par
Assertion (2). 
In this case one has $(d f_{E_0})_{q_0} = -1$ and $(d f_{E_1})_{q_0} = - \delta$, 
because $f_{E_0}$ fixes $q_0$ and exchanges $q_{\pm}$. 
As in the last paragraph one gets $(d f_{E_j})_{q_j} = -\delta^{-j}$, 
$(d f_{E_{j+1}})_{q_j} = - \delta^{j+1}$ successively for $j = 1, \dots, k-1$, until arriving 
at a unique fixed point $q_k \in E_k$ different from $q_{k-1}$, at which  
$(d f)_{q_k}$ has eigenvalues $-\delta^{-k}$ and $-\delta^{k+1}$ (see Figure \ref{fig:excep}). 
Do not forget that $f_{E_0}$ has one more fixed point $q_0' \in E_0$ different from 
$q_0$, at which $(d f)_{q_0'}$ has eigenvalues $-1$ and $-\delta$.    
The $k+2$ points $q_0', q_0, q_1, \dots, q_k$, which are transverse, 
are all of the isolated fixed points on $E_0 \cup E_1 \cup \cdots \cup E_k$ 
and the sum of their $\nu$-indices is given by   
$$
\varLambda_k^- 
:= \nu_{q_0'}(f) + \sum_{j=0}^k \nu_{q_j}(f) 
= \dfrac{1}{2(1+\delta)} + \sum_{j=0}^k \dfrac{1}{(1+\delta^{-j})(1+\delta^{j+1})} 
= \dfrac{1}{2(1+\delta)} + \dfrac{1+ \delta + \cdots + \delta^k}{2(1+\delta^{k+1})}.
$$
If $\varGamma$ is of type $\rD_n$ then $\mu(f, \cE') = (n-3) + 2 = n-1$ and 
$\nu(f, \cE') = \varLambda_{n-3}^-$, since contributing to $\nu(f, \cE')$ is only the 
longest arm of $\varGamma$, with length $n-3$; this  yields \eqref{eqn:numD}. 
If $\varGamma$ is of type $\rE_6$ then $\mu(f, \cE') = 1+2 = 3$ and 
$\nu(f, \cE') = \varLambda_1^+$, since contributing to $\nu(\cE')$ is only the shortest 
arm of $\varGamma$, with length $1$; this yields \eqref{eqn:numE}. 
\hfill $\Box$ \par\medskip
Some calculations related to Lemma \ref{lem:DE} are made in \cite[\S 9.3]{IT} 
for a couple of examples. 
Now Lemma \ref{lem:DE} provides a thorough result of this sort in a 
unified manner. 
When $\cE'$ is a connected component of Dynkin type $\rA$, things are much 
subtler due to the possible occurrence of a {\sl multiple} isolated fixed point 
on $\cE'$. 
In this article no attempt is made to develop a general theory for components 
of type $\rA$.   
Instead, a particular case is discussed in \S \ref{sec:P2}, where only one  
$\rA_1$-component is present. 
Even in this case situations are already hard and interesting.  
The results in the next section will be needed in this context, 
although they are important in their own light.                
\section{Grothendieck Residues} \label{sec:res} 
Evaluating the residue $\nu_p(f)$ in \eqref{eqn:g-res} at a multiple fixed point $p$ 
is usually difficult, but when $p$ lies on an invariant curve there are cases 
where this task is more tractable.   
The aim of this section is to discuss such situations, or more precisely, to do so 
in a dynamical context, not only with the map $f$ alone but also with its iterates.      
\begin{definition} \label{def:Fixe}
A fixed point $p \in X$ of $f$ is said to be {\sl exceptional} if 
$2 \le \mu_p(f) < \infty$ and there exists an $f$-invariant curve $E$  
passing through $p$. 
For such a curve $E$, since $\mu_p(f) \ge 2$ and 
$\det (d f)_p = \delta(f)$, one has either (i) $(d f_E)_p = 1$; or 
(ii) $(d f_E)_p = \delta(f)$, where $f_E :=f|_E$ is the M\"{o}bius transformation 
$f$ induces on $E \cong \bP^1$.  
Let $\Fixe(f)$ be the set of all exceptional fixed points of $f$.  
Condition $\mu_p(f) < \infty$ implies $\Fixe(f) \subset \Fixi(f)$.      
We say that $p \in \Fixe(f)$ is {\sl of type I} if $p$ admits a curve 
$E$ of type (i); and {\sl of type I\!I} if $p$ admits no curve of type (i) 
but a curve $E$ of type (ii). 
Let $\Fixe(f) = \Fixe_{\rI}(f) \amalg \Fixe_{\rII}(f)$ be the decomposition  
according to the types.    
\end{definition} 
\begin{lemma} \label{lem:Fixe} 
For any integer $n \ge 1$ let $f^n := f \circ \cdots \circ f$ be the 
$n$-th iterate of $f$. 
Then  
\begin{equation} \label{eqn:Fixe}
\Fixe_{\rI}(f) \subset \Fixe_{\rI}(f^n), \qquad 
\Fixe_{\rII}(f) \subset \Fixe_{\rII}(f^n) \qquad \mbox{for any} \quad n \ge 1.  
\end{equation}
\end{lemma}
{\it Proof}. 
First let $p \in \Fixe_{\rI}(f)$ and $E$ be an $f$-invariant but not fixed curve 
passing through $p$ with $(d f_E)_p = 1$. 
The M\"{o}bius transformation $f_E$ and its iterates $f^n_E$ can then 
be expressed as  
\begin{equation} \label{eqn:parabolic}
f_E(z_1) = \dfrac{z_1}{1+ z_1} \qquad \mbox{and hence} \qquad 
f^n_E(z_1) = \dfrac{z_1}{1+ n z_1},    
\end{equation}
in terms of a suitable coordinate $z_1$ on $E$ such that $z_1 = 0$ at $p$.   
Thus $E$ is not a fixed curve of $f^n$ for any $n \ge 1$.  
Suppose that for some $n \ge 2$ there is a fixed curve $C$ of $f^n$ 
passing through $p$.   
By Lemma \ref{lem:transv}, $C$ intersects $E$ transversally in $p$, so 
$(d f_E^n)_p = 1$ and $\det (d f^n)_p = \delta^n$ imply 
$1 = (d f^n_C)_p = \delta^n$, which contradicts the assumption that $\delta$ 
is conjugate to a Salem number. 
Therefore $p \in \Fixi(f^n)$ and so $p \in \Fixe_{\rI}(f^n)$ for any $n \ge 1$.     
\par
Next let $p \in \Fixe_{\rII}(f)$ and $E$ be an $f$-invariant but not fixed curve 
passing through $p$ such that $(d f_E)_p = \delta$. 
Since $\delta$ is not a root of unity, $E$ is not a fixed curve of $f^n$ 
for any $n \ge 1$.  
Suppose that for some $n \ge 2$ there is an $f^n$-invariant curve 
$C \neq E$ passing through $p$. 
By Lemma \ref{lem:transv}, $C$ meets $E$ transversally in $p$, hence  
$(d f_C^n)_p = 1$ by $(d f_E^n) = \delta^n = \det(d f^n)_p$. 
Similarly, if $C$ is $f$-invariant then $(d f_C)_p = 1$ and   
$p \in \Fixe_{\rI}(f)$, contradicting $p \in \Fixe_{\rII}(f)$. 
So $C' := f(C)$ is a different $f^n$-invariant curve with 
$(d f^n_{C'})_p = \delta^n$, which intersects $C$ transversally in $p$ 
by Lemma \ref{lem:transv}. 
As $T_pX = T_p C \oplus T_pC'$, we have $\delta^n = \det(d f^n)_p = 
\delta^n \cdot \delta^n$, i.e. $\delta^n = 1$, again a contradiction. 
Thus $f^n$ admits no invariant curve through $p$ 
other than $E$ and hence $p \in \Fixe_{\rII}(f^n)$ for any $n \ge 1$.          
\hfill $\Box$
\subsection{Exceptional Fixed Points of Type I} \label{ss:2nd} 
The case of type I is simpler to deal with than that of type I\!I, 
so we begin with the former.  
\begin{theorem} \label{thm:mu}
If $p \in \Fixe_{\rI}(f)$ then we have $p \in \Fixe_{\rI}(f^n)$ and 
\begin{equation} \label{eqn:mu} 
\mu(f^n) = 2, \qquad 
\nu_p(f^n) = \dfrac{1 + \delta^n}{(1-\delta^n)^2} \qquad 
\mbox{for any} \quad n \ge 1.     
\end{equation}
\end{theorem} 
\par
This theorem will be established after Lemma \ref{lem:mu}.  
Formulas \eqref{eqn:mu_C} and \eqref{eqn:mu} tell us that 
an irreducible fixed curve and an exceptional fixed point of type I  
have the same holomorphic index. 
\par
Suppose $p \in \Fixe_{\rI}(f)$ and let $E$ be the curve in Definition \ref{def:Fixe}. 
Along with $(d f_E)_p = 1$ the tangent map $(d f)_p$ has eigenvalue $\delta$. 
Take a local chart $z = (z_1, z_2)$ around $p = (0, 0)$ such that $z_2$ is in the 
eigen-direction of eigenvalue $\delta$, $E = \{ z_2 = 0\}$ and $z_1$ is a 
coordinate on $E$ such that $f_E$ is normalized as in \eqref{eqn:parabolic}.   
Let $(f_1, f_2)$ be the local representation for $f$ in this chart.  
Since $f$ preserves $E$, there exist $g_1(z)$, $g_2(z) \in \bC\{ z \}$ 
such that  
\begin{equation} \label{eqn:f1}
f_1(z) = \dfrac{z_1}{1+z_1}  + z_2 \, g_1(z), \qquad 
f_2(z) = z_2 \{ \delta + g_2(z) \}, 
\qquad g_1(0, 0) = g_2(0, 0) = 0.  
\end{equation}
\begin{lemma} \label{lem:mult1} 
For any $p \in \Fixe_{\rI}(f)$ we have 
$\frak{a} = (z_1^2, z_2)$ for the ideal in \eqref{eqn:index} 
and hence $\mu_p(f) = 2$. 
\end{lemma}
{\it Proof}. 
First we have $z_2 \in \frak{a}$, since 
$z_2 - f_2(z) = z_2 \{ 1-\delta- g_2(z) \} \in \frak{a}$ and 
$1-\delta- g_2(z) \in \bC\{z\}^{\times}$. 
Secondly we have $z_1^2 \in \frak{a}$, since 
$z_1 - f_1(z)  + z_2 \, g_1(z) = z_1^2/(1+z_1) \in \frak{a}$ and 
$1/(1+z_1) \in \bC\{z_1\}^{\times}$. 
Thus $(z_1^2, z_2) \subset \frak{a}$ and the converse inclusion is obvious. 
As $1$ and $z_1$ form a basis of $\bC\{ z \}/\frak{a}$, we have $\mu_p(f) = 2$.  
\hfill $\Box$ \par\medskip 
Let $\ve \approx 0$ be a small parameter and set $t := 1/(1-\ve)$. 
Consider a perturbation $f^{\ve}$ of $f$ defined by 
$$
f^{\ve}(z) := f(t z_1, t^{-1} z_2), \qquad (\ve, z) \approx (0, 0, 0). 
$$
\begin{lemma} \label{lem:p-fixed1} 
For every $\ve \approx 0$ with $\ve \neq 0$ tha map $f^{\ve}$ has exactly two 
fixed points $p = (0, 0)$ and $q^{\ve} := (\ve, 0)$ in a small neighborhood of 
$(\ve, z_1, z_2) = (0, 0, 0)$. 
Moreover, $\delta^{\ve} := \det (d f^{\ve})_{q^{\ve}} = \delta(1 + a \ve^2)$ 
with $a = O(1)$ as $\ve \to 0$.  
\end{lemma}
{\it Proof}. 
In view of \eqref{eqn:f1} the equation $z_2 - f_2^{\ve}(z) = 0$ reads 
$z_2\{1 - t \delta - t g_2(t z_1, t^{-1} z_2)\} = 0$, which yields $z_2 = 0$ as  
$1 - t \delta - t g_2(t z_1, t^{-1} z_2) \approx 1 - \delta \neq 0$ for  
$(\ve, z_1, z_2) \approx (0, 0, 0)$.  
Again by \eqref{eqn:f1} putting $z_2 = 0$ into $z_1 - f_1^{\ve}(z) = 0$ 
gives $z_1(z_1-\ve)/(1-\ve)=0$, i.e. $z_1 = 0$, $\ve$. 
Thus the fixed points of $f^{\ve}$ are exactly $p = (0, 0)$ and $q^{\ve} = (\ve, 0)$. 
\par
If $\eta$ is a nowhere vanishing holomorphic $2$-form on $X$, then  
equation $f^* \eta = \delta \cdot \eta$ is represented as 
\begin{equation} \label{eqn:J}
J_f(z) = \delta \, \dfrac{h(z)}{h(f(z))} \qquad \mbox{with} \qquad 
\eta = h(z)\, d z_1 \wedge d z_2,  
\end{equation}
where $J_f(z)$ is the Jacobian of $f$. 
Substituting $z = (t \ve, 0)$ into \eqref{eqn:J} and using $f(t \ve, 0) = (\ve, 0)$, 
we have  
$$
\delta^{\ve} := \det (d f^{\ve})_{q^{\ve}} = J_f(t \ve, 0) = \delta \, 
\dfrac{h(t \ve, 0)}{h(f(t \ve, 0))} = \delta \, 
\dfrac{h(\ve (1-\ve)^{-1}, 0)}{h(\ve, 0)} = \delta \left\{ 1 + O(\ve^2) \right\}, 
$$
since $h(\ve, 0) = c_0 \{ 1 + c_1 \ve + O(\ve^2)\}$ for some constants 
$c_0 \in \bC^{\times}$ and $c_1 \in \bC$.  \hfill $\Box$
\begin{lemma} \label{lem:mu} 
For any $p \in \Fixe_{\rI}(f)$ formula \eqref{eqn:mu} holds for $n = 1$. 
\end{lemma}
{\it Proof}. 
It is easy to see that $(d f^{\ve})_p$ has eigenvalues $(d f^{\ve}_E)_p = t$ 
and $\delta t^{-1}$. 
One then has $(d f^{\ve}_E)_{q^{\ve}} = t^{-1}$, since $q^{\ve} \in E$ is the 
other fixed point of the M\"{o}bius transformation $f^{\ve}_E$. 
So $(d f^{\ve})_{q^{\ve}}$ has eigenvalues $t^{-1}$ and $t \delta^{\ve}$.  
Let $\omega^{\ve}$ be the $2$-form in \eqref{eqn:g-res} for $f^{\ve}$. 
By continuity principle \cite[\S 5.1]{GH} the residue $\nu_p(f)$ 
is given as the limit of  
\begin{align*}
\Res_p \, \omega^{\ve} + \Res_{q^{\ve}} \, \omega^{\ve} 
&= \dfrac{1}{(1-t)(1-t^{-1} \delta)} + \dfrac{1}{(1-t^{-1})(1-t \delta^{\ve})} 
\\[2mm]
&= \dfrac{ (1-\ve)(1+\delta+ a \delta \ve) }{(1-\delta+ \delta \ve) 
(1-\delta -\ve - a \delta \ve^2)} 
\to \dfrac{1+ \delta }{(1-\delta)^2} 
\qquad \mbox{as} \quad \ve \to 0,   
\end{align*}
where formula \eqref{eqn:h-index} and Lemma \ref{lem:p-fixed1} are 
used in the first and second equalities respectively. 
\hfill $\Box$ \par\medskip
{\it Proof of Theorem $\ref{thm:mu}$}. 
It is an immediate consequence of 
Lemmas \ref{lem:Fixe}, \ref{lem:mult1} and \ref{lem:mu}.  \hfill $\Box$
\subsection{Exceptional Fixed Points of Type I\!I} \label{ss:1st} 
Let $p \in X$ be a fixed point of $f$ lying on an invariant curve $E$ such that  
$(d f_E)_p = \delta$.  
Note that $f_E$ has a fixed point $q \in E$ different from $p$, 
at which $(d f_E)_q = \delta^{-1}$.     
Along with $\delta$ the tangent map $(d f)_p$ has eigenvalue $1$. 
Take a local chart $z = (z_1, z_2)$ around $p = (0, 0)$ such that 
$z_1$ is in the eigen-direction of eigenvalue $1$, $E = \{ z_1 = 0\}$, and 
$z_2$ gives a coordinate on $E \setminus \{ q \} \cong \bC$ 
with $q$ located at $z_2 = \infty$, so that $f_E$ is normalized as 
$f_E(z_2) = \delta z_2$.  
Let $(f_1, f_2)$ be the local representation for $f$ in this chart.  
Since $f$ preserves $E$, we can write 
\begin{equation} \label{eqn:f2}
f_1(z) = z_1 \{ 1 + g_1(z) \}, \qquad f_2(z) = \delta \{ z_2 + z_1 \, g_2(z) \}, 
\qquad g_1(0, 0) = g_2(0, 0) = 0, 
\end{equation}
for some $g_1(z)$, $g_2(z) \in \bC\{ z \}$.  
Provide $z_1$ and $z_2$ with orders $1$ and $2$ respectively. 
Put 
\begin{equation} \label{eqn:g2}
g_1(z) = \sum_{i, j = 0}^{\infty} a_{i j} \, z_1^i z_2^j, \qquad 
g_2(z) = \sum_{i, j = 0}^{\infty} b_{i j} \, z_1^i z_2^j, \qquad a_{00} = b_{00} = 0. 
\end{equation}
\begin{lemma} \label{lem:mult2} 
We have $\mu_p(f) = 2$ if and only if $a_{10} \neq 0$, in which case 
$\frak{a} = (z_1^2, z_2)$ for the ideal in \eqref{eqn:index}.     
\end{lemma}
{\it Proof}. 
First we show that $a_{10} \neq 0$ implies $\frak{a} = (z_1^2, z_2)$ and 
$\mu_p(f) = 2$. 
In what follows $u_j(z_1)$ and $v_j(z)$ stand for various elements in 
$\bC\{ z_1\}$ and $\bC\{ z \}$ respectively.  
Observe that $z_1 - f_1(z) = z_1 \{ z_1 u_1(z_1) + z_2 v_1(z) \} \in \frak{a}$ 
with $u_1(z_1) \in \bC\{ z_1\}^{\times}$.  
Multiplying it by $u_1(z_1)^{-1}$ yields  
$v_2(z) := z_1 \{ z_1 + z_2 v_3(z) \} \in \frak{a}$ with 
$v_3(z) := u_1(z_1)^{-1} v_1(z)$. 
One also has $z_2 - f_2(z) = z_2 \{ 1-\delta + z_1 v_4(z) \} + 
z_1^{k+2} u_2(z_1) \in \frak{a}$ for some $k \ge 0$. 
So $z_2 - f_2(z) - z_1^k u_2(z_1) \, v_2(z) = z_2 v_4(z) \in \frak{a}$ 
with $v_4(z) = 1-\delta + z_1 \{ v_4(z) - z_1^k u_2(z_1) \, v_3(z) \} 
\in \bC\{ z \}^{\times}$. 
Thus $z_2 \in \frak{a}$ and $z_1^2 = v_2(z) - z_1 z_2 \, v_3(z) \in \frak{a}$, 
hence $\frak{a} = (z_1^2, z_2)$. 
Therefore $1$ and $z_1$ give a basis of $\bC\{ z \}/\frak{a}$, so we 
have $\mu_p(f) = 2$.      
\par
Next we show that $a_{10} = 0$ implies $\mu_p(f) \ge 3$. 
It suffices to prove the following two claims.   
\begin{enumerate}
\setlength{\itemsep}{-1pt}
\item If $b_{10} \neq 0$ then $1$, $z_1$, $z_2$ are linearly independent 
in $\bC\{ z \}/\frak{a}$. 
\item If $b_{10} = 0$ then $1$, $z_1$, $z_1^2$ are linearly independent 
in $\bC\{ z \}/\frak{a}$. 
\end{enumerate}
Note that $\ord (z_1 - f_1(z)) \ge 3$ and $\ord (z_2 - f_2(z)) \ge 2$. 
To show claim (1), suppose that 
$\alpha_0 + \alpha_1 z_1 + \alpha_2 z_2 \in \frak{a}$ with 
$\alpha_0$, $\alpha_1$, $\alpha_2 \in \bC$. 
Since any element of $\frak{a}$ has order at least $2$, we have 
$\alpha_0 = \alpha_1 = 0$ and $\alpha_2 z_2 \in \frak{a}$.  
Thus $\alpha_2 z_2 = u(z) \cdot (z_1 - f_1(z)) + v(z) \cdot (z_2 - f_2(z))$ for 
some $u(z)$, $v(z) \in \bC\{ z \}$. 
Its order $2$ component yields $\alpha_2 z_2 = v(0, 0) 
\{(1-\delta) z_2 - \delta b_{10} z_1^2 \}$. 
As $b_{10} \neq 0$ we have $v(0, 0) = 0$ and hence $\alpha_2 = 0$. 
To show claim (2), suppose that 
$\alpha_0 + \alpha_1 z_1 + \alpha_2 z_1^2 \in \frak{a}$ 
with $\alpha_0$, $\alpha_1$, $\alpha_2 \in \bC$. 
Since $\ord (z_1 - f_1(z))|_{z_2=0} \ge 3$ and $\ord (z_2 - f_2(z))|_{z_2=0} \ge 3$, 
any element of $\frak{a}|_{z_2 = 0}$ has order at least $3$.  
This implies $\alpha_0 = \alpha_1 = \alpha_2 = 0$. \hfill $\Box$ 
\begin{theorem} \label{thm:iterate} 
If $p \in \Fixe_{\rII}(f)$ and $\mu_p(f) =2$ then we have    
$p \in \Fixe_{\rII}(f^n)$ and 
\begin{equation} \label{eqn:iterate}
\mu_p(f^n) = 2, \qquad  
\nu_p(f^n) = 
\dfrac{n-1 + (n+1) \delta^n + 
\left(1+ \delta + \cdots + \delta^{n-1} \right) \theta}{ n \, (1-\delta^n)^2} 
\qquad \mbox{for any} \quad n \ge 1,  
\end{equation}
where in terms of some leading coefficients in \eqref{eqn:g2} the quantity 
$\theta$ is defined by 
\begin{equation} \label{eqn:delta}
\theta :=  \dfrac{(1-\delta) a_{20} + \delta \, a_{01} b_{10}}{ (a_{10})^2  }. 
\end{equation}
\end{theorem}
\par
We establish this theorem by providing four lemmas, where we work with 
$f$ in Lemmas \ref{lem:b01}--\ref{lem:res2} and proceed to its iterates $f^n$ 
in Lemma \ref{lem:iterate2}. 
Hereafter we assume \eqref{eqn:f2}, \eqref{eqn:g2} and $\mu_p(f) = 2$ 
without further comment.      
Rescaling $z \to (\lambda z_1, z_2)$ with $\lambda \in \bC^{\times}$ takes 
$f_1(z) \to \lambda^{-1} f_1(\lambda z_1, z_2)$ and $f_2(z) \to f_2(\lambda z_1, z_2)$, 
hence induces the change of coefficients $a_{i j} \to \lambda^i a_{i j}$ and 
$b_{i j} \to \lambda^{i+1} b_{i j}$. 
We can take $\lambda = (a_{10})^{-1}$ to get a normalization  
\begin{equation} \label{eqn:norm}
a_{10} = 1. 
\end{equation}  
\begin{lemma} \label{lem:b01} 
We have $b_{01} = -2$ under the normalization \eqref{eqn:norm}.  
\end{lemma}
{\it Proof}. 
This follows from equation \eqref{eqn:J}. 
Indeed, substituting \eqref{eqn:f2} and \eqref{eqn:g2} into it,  
we observe   
$$
\mbox{LHS of \eqref{eqn:J}} = \delta \{ 1 + (b_{01} + 2) z_1 + O_2 \}, 
\qquad \mbox{RHS of \eqref{eqn:J}} = \delta \{ 1 + O_2 \},  
$$
where $O_2$ stands for various terms of oder at least $2$. 
Comparing the first order terms yields $b_{01} + 2 = 0$. 
\hfill $\Box$ \par\medskip
Let $\ve \approx 0$ be a small parameter and set $t := 1-\ve$. 
Consider a perturbation $f^{\ve}$ of $f$ defined by 
$$
f^{\ve}(z) := f(t z_1, t^{-1} z_2), \qquad (\ve, z) \approx (0, 0, 0). 
$$
It is obvious that the origin $p = (0, 0)$ is a fixed point of $f^{\ve}$. 
Let us find another fixed point. 
\begin{lemma} \label{lem:p-fixed2} 
For every $\ve \approx 0$, $\ve \neq 0$, the map $f^{\ve}$ has exactly two 
fixed points $p = (0, 0)$ and $q^{\ve} = (w_1(\ve), w_2(\ve))$ in a small 
neighborhood of  $(\ve, z_1, z_2) = (0, 0, 0)$.  
Under \eqref{eqn:norm} the coordinates of $q^{\ve}$ admit an expansion    
\begin{subequations} \label{eqn:p-fixed2} 
\begin{alignat}{2} 
w_1(\ve) &= \ve + A_2 \ve^2 + O(\ve^3), \qquad & 
A_2 &= 2 - a_{20} - \dfrac{\delta \, a_{01} b_{10}}{1-\delta}, \label{eqn:wA}  \\[1mm]
w_2(\ve) &= B_2 \ve^2 + O(\ve^3),  \qquad & B_2 &= \dfrac{\delta \, b_{10}}{1-\delta}.  
\label{eqn:wB} 
\end{alignat} 
\end{subequations}
\end{lemma}
{\it Proof}. 
Put $F_1^{\ve}(z) := t \, g_1(t z_1, t^{-1} z_2) - \ve$ and 
$F_2^{\ve}(z) := (\delta t^{-1} -1) z_2 + \delta \, t z_1 \, g_2(t z_1, t^{-1} z_2)$.  
Then $f^{\ve}(z) = z$ is equivalent to 
$z_1 F^{\ve}_1(z) = F^{\ve}_2(z) = 0$. 
If $z_1 = 0$ then $(\delta t^{-1} - 1) z_2 = 0$ and hence $z_2 = 0$, 
as $\delta t^{-1} - 1 \approx \delta-1 \neq 0$ for $\ve \approx 0$.  
Thus any fixed point other than $p$ is a solution to the equations  
$F^{\ve}_1(z) = F^{\ve}_2(z) = 0$. 
Since 
$$
\dfrac{\partial F^0}{\partial z}(0, 0) = 
\begin{pmatrix} 1 & a_{01} \\[1mm] 0 & \delta-1\end{pmatrix},  \qquad 
\left. \dfrac{\partial F^{\ve}}{\partial \ve}(0, 0) \right|_{\ve = 0} = - 
\begin{pmatrix} 1 \\[1mm] 0 \end{pmatrix} \qquad 
\mbox{under \eqref{eqn:norm}},  
$$ 
the implicit function theorem implies that near $(\ve, z_1, z_2) = (0, 0, 0)$ 
there exists a unique solution $(z_1, z_2) = (w_1(\ve), w_2(\ve))$ such that 
$w_1(0) = w_2(0) = 0$. 
It satisfies $w_1'(0) = 1$, $w_2'(0) = 0$, so admits an expansion as in 
\eqref{eqn:p-fixed2} for some constants $A_2$ and $B_2$.  
Using \eqref{eqn:f2}, \eqref{eqn:g2}, \eqref{eqn:norm} and Lemma \ref{lem:b01}, 
we have $F_1^{\ve}(w_1(\ve), w_2(\ve)) 
= (a_{20}-2 + A_2 + a_{01} B_2) \ve^2 + O(\ve^3) = 0$ and    
$F_2^{\ve}(w_1(\ve), w_2(\ve)) = 
\{ (\delta-1) B_2 + \delta b_{10} \} \ve^2 + O(\ve^3) = 0$, hence  
$a_{20}-2 + A_2 + a_{01} B_2 = 0 = (\delta-1) B_2 + \delta b_{10}$. 
This determines $A_2$ and $B_2$ as in \eqref{eqn:wA} and \eqref{eqn:wB} 
respectively. \hfill $\Box$
\begin{lemma} \label{lem:res2} In terms of $\theta$ in \eqref{eqn:delta} 
the holomorphic local index in \eqref{eqn:g-res} can be expressed as   
\begin{equation} \label{eqn:res1} 
\nu_p(f) = \dfrac{2 \delta + \theta}{(1-\delta)^2}.  
\end{equation}
\end{lemma}
{\it Proof}.  We keep the normalization \eqref{eqn:norm}. 
Since $(d f^{\ve})_p$ has eigenvalues $t$ and $\delta t^{-1}$,   
$$
\det \left( I - (d f^{\ve})_p \right) = 
(1-t)(1- \delta t^{-1}) = \dfrac{\ve(1-\delta-\ve) }{1-\ve}. 
$$
Calculation of $\det \left( I - (d f^{\ve})_{q^{\ve}} \right)$ is much 
harder, but Lemmas \ref{lem:b01} and \ref{lem:p-fixed2} together 
with \eqref{eqn:f2}, \eqref{eqn:g2}, \eqref{eqn:norm} yield    
$$
\det \left( I - (d f^{\ve})_{q^{\ve}} \right) = 
\ve \left[ \delta-1 -\{ \delta +(1-\delta) a_{20} + \delta \, a_{01} b_{10} \} \ve 
+ O(\ve^2) \right]. 
$$
Let $\omega^{\ve}$ be the $2$-form in \eqref{eqn:g-res} for $f^{\ve}$. 
By continuity principle \cite[\S 5.1]{GH} the residue $\nu_p(f)$ 
is given as the limit of  
\begin{align*}
\Res_p \, \omega^{\ve} + \Res_{q^{\ve}} \, \omega^{\ve} 
&= \dfrac{1}{\det \left( I - (d f^{\ve})_p \right)} + 
\dfrac{1}{ \det \left( I - (d f^{\ve})_{q^{\ve}} \right) } \\[2mm] 
&= \dfrac{ 2 \delta + (1-\delta) a_{20} + \delta \, a_{01} b_{10} -  
\{ \delta +(1-\delta) a_{20} + \delta \, a_{01} b_{10} \} \ve + 
O(\ve^2) }{(1- \delta - \ve)[1 - \delta + \{\delta + (1-\delta) a_{20} + 
\delta \, a_{01} b_{10} \} \ve + O(\ve^2)]} \\[2mm]
&\to \dfrac{2 \delta + \theta}{(1-\delta)^2} 
\qquad \mbox{as} \quad \ve \to 0 \quad \mbox{with} \quad 
\theta := (1-\delta) a_{20} + \delta \, a_{01} b_{10}. 
\end{align*}
Removing the normalization \eqref{eqn:norm} we obtain  
formula \eqref{eqn:res1} with $\theta$ defined in \eqref{eqn:delta}. 
\hfill $\Box$ \par\medskip
For any $n \ge 1$ the $n$-th iterate $f^n$ can also be represented in the form 
\eqref{eqn:f2}-\eqref{eqn:g2} upon replacing $\delta$ by $\delta^n$ in \eqref{eqn:f2} 
and rewriting the coefficients $a_{ij}$ and $b_{ij}$ as $a_{ij}^{(n)}$ and $b_{ij}^{(n)}$ in 
\eqref{eqn:g2}, respectively.   
\begin{lemma} \label{lem:iterate2} 
For any $n \ge 1$ one has $\mu_p(f^n) = 2$ and hence $p$ is an isolated 
fixed point of $f^n$.    
Moreover,  
\begin{equation} \label{eqn:n-delta}
\theta^{(n)} := \dfrac{ (1-\delta^n) a_{20}^{(n)} + 
\delta^n \, a_{01}^{(n)} b_{10}^{(n)} }{ (a_{10}^{(n)})^2} = 
\dfrac{(1-\delta^n) \{(n-1)(1-\delta) + \theta \}}{n(1-\delta)} \qquad 
\mbox{for any} \quad n \ge 1.  
\end{equation}
\end{lemma}
{\it Proof}. 
The obvious composition rule $f^{n+1} = f \circ f^n$ then leads to a system of 
recurrence relations  
\begin{alignat*}{2}
a_{10}^{(n+1)} &= a_{10}^{(n)} + 1, \qquad & 
a_{01}^{(n+1)} &= a_{01}^{(n)} + \delta^n a_{01},  \\
b_{10}^{(n+1)} &= b_{10}^{(n)} + \delta^{-n} b_{10}, \qquad & 
a_{20}^{(n+1)} &= a_{20}^{(n)} + a_{20} + 2 a_{10}^{(n)} + \delta^n a_{01} b_{10}^{(n)},  
\end{alignat*}
where the normalization \eqref{eqn:norm} is employed. 
This system is readily settled as   
\begin{alignat*}{2}
a_{10}^{(n)} &= n \ge 1, \qquad & 
a_{01}^{(n)} &= \dfrac{(1-\delta^n) a_{01}}{1-\delta},  \\
b_{10}^{(n)} &= \dfrac{(1-\delta^n) b_{10}}{\delta^{n-1} (1-\delta)}, \qquad & 
a_{20}^{(n)} &= n(n-1) + n \, a_{20} + \dfrac{\delta}{1-\delta} \left\{ 
n-1 - \dfrac{\delta(1-\delta^{n-1})}{1-\delta} \right\} a_{01} b_{10}. 
\end{alignat*}
Lemma \ref{lem:mult2} shows that $\mu_p(f^n) = 2 < \infty$ and hence 
$p$ is an isolated fixed point of $f^n$ for every $n \ge 1$.   
Substituting the above data into the definition of $\theta^{(n)}$, we find that 
a fine cancellation occurs to yield \eqref{eqn:n-delta}.  \hfill $\Box$ \par\medskip
{\it Proof of Theorem $\ref{thm:iterate}$}. 
It is clear from Lemma \ref{lem:iterate2} that $p$ is an exceptional fixed point 
of type I\!I relative to $E$. 
Lemma \ref{lem:res2} imples $\nu_p(f^n) = 
(2 \delta^n + \theta^{(n)})/(1-\delta^n)^2$, which combined with 
\eqref{eqn:n-delta} yields formula \eqref{eqn:iterate}. \hfill $\Box$   
\section{Siegel Disks} \label{sec:SD}
Let $f : X \to X$ be a K3 surface automorphism satisfying the conditions 
(C1) and (C2) at the beginning of \S \ref{sec:fpf}. 
In \cite[Proposition 9.1]{IT} we give a criterion for a given fixed point of $f$ to be 
the center of a Siegel disk or to be a hyperbolic fixed point.  
For later use we have to extend it a little bit. 
Let $p \in X$ be a fixed point of $f$.
Then the eigenvalues of the tangent map $(d f)_p : T_pX \to T_pX$ can be 
represented as 
\begin{equation} \label{eqn:ev1}
\alpha_1 := \delta^{\frac{1}{2}} \alpha, \qquad 
\alpha_2 := \delta^{\frac{1}{2}} \alpha^{-1}  \qquad 
\mbox{for some} \quad  \alpha \in \bC^{\times}, 
\end{equation}
where the branch of $\delta^{\frac{1}{2}}$ is specified by 
$\Re(\delta^{\frac{1}{2}}) > 0$ for the sake of definiteness. 
Let $\tau$ be the special trace of $f$.  
\begin{lemma} \label{lem:mi1} 
Suppose that there exists a rational functions $P(w) \in \bQ(w)$ such that  
$(\alpha + \alpha^{-1})^2 = P(\tau)$. 
\begin{enumerate}
\setlength{\itemsep}{-1pt}
\item If $0 \le P(\tau) \le 4$ then $p$ is the center of a Siegel disk, provided 
either $(\mathrm{i})$ $\tau$ admits a conjugate $\tau'$ such that 
$-2 < \tau' < 2$ and $P(\tau') > 4;$ or $(\mathrm{ii})$ $P(\tau)$ is 
not an algebraic integer.  
\item If $P(\tau) > 4$ then $p$ is a hyperbolic fixed point. 
\end{enumerate} 
\end{lemma}
{\it Proof}. 
The cases of (1)-(i) and (2) are proved in \cite[Proposition 9.1]{IT}.  
Under the condition that $\alpha_1$, $\alpha_2 \in \overline{\bQ} \cap S^1$,  
the fixed point $p$ is the center of a Siegel disk if and only if 
$\alpha_1$ and $\alpha_2$ are multiplicatively independent (MI) 
(see \cite[Theorem 5.1]{McMullen1}). 
In case (1)-(ii) it follows from $(\alpha + \alpha^{-1})^2 = P(\tau)$ and 
$0 \le P(\tau) \le 4$ that $\alpha_1$, $\alpha_2 \in \overline{\bQ} \cap S^1$, 
while assumption (ii) implies that $\alpha$ is {\sl not} an algebraic unit.  
Suppose that $\alpha_1^m \alpha_2^n = \delta^{\frac{1}{2}(m+n)} \alpha^{m-n} =1$, 
that is, $\alpha^{n-m} = \delta^{\frac{1}{2}(m+n)}$ for 
some $m$, $n \in \bZ$.  
If $n-m \neq 0$ then $\alpha$ must be an algebraic unit, since so is $\delta$, 
but this is impossible.  
Thus one has $n-m = 0$ and $\delta^{\frac{1}{2}(m+n)} = 1$, but the latter equation 
yields $m+n = 0$, because $\delta$ is not a root of unity as a conjugate to 
a Salem number.  
Thus $m = n = 0$ and hence $\alpha_1$ and $\alpha_2$ are MI.   
\hfill $\Box$ \par\medskip
It is sometimes more convenient to express the eigenvalues of $(d f)_p$ in the form  
\begin{equation} \label{eqn:ev2} 
\beta_1 := \beta, \qquad  \beta_2 := \delta \beta^{-1} \qquad 
\mbox{for some} \quad \beta \in \bC^{\times}.  
\end{equation}
An obvious variant of Lemma \ref{lem:mi1} in this situation is the following lemma, 
whose proof is safely omitted. 
\begin{lemma} \label{lem:mi2} 
Suppose that there exists a rational functions $Q(w) \in \bQ(w)$ such that  
$\beta + \beta^{-1} = Q(\tau)$.  
\begin{enumerate}
\setlength{\itemsep}{-1pt}
\item If $|Q(\tau)| \le 2$ then $p$ is the center of a Siegel disk, 
provided either $(\mathrm{i})$ $\tau$ admits a conjugate $\tau'$ such that 
$-2 < \tau' < 2$ and $|Q(\tau')| > 2;$ or $(\mathrm{ii})$ $Q(\tau)$ is not 
an algebraic integer. 
\item If $|Q(\tau)| > 2$ then $p$ is a hyperbolic fixed point. 
\end{enumerate} 
\end{lemma}
\begin{remark} \label{rem:mi} 
In Lemmas $\ref{lem:mi1}$ and $\ref{lem:mi2}$, suppose that $p$ lies on a 
$(-2)$-curve $E \subset \cE$.    
Then $f_E$ admits a unique fixed point $p' \in E$ other than $p$. 
If $p$ is the center of a Siegel disk, then so is $p'$.  
If $p$ is a hyperbolic fixed point, then so is $p'$. 
Indeed, if the eigenvalues of $(d f)_p$ are given by \eqref{eqn:ev1} with 
$(d f_E)_p = \alpha_1$, then those of $(d f)_{p'}$ are 
$\alpha_1' := \alpha_1^{-1} = \delta^{-\frac{1}{2}} \alpha^{-1}$ and 
$\alpha_2' := \delta \alpha_1 = \delta^{\frac{3}{2}} \alpha$. 
Note that $\alpha_1$ and $\alpha_2$ are MI, if and only if $\alpha$ and 
$\delta^{\frac{1}{2}}$ are MI, if and only if $\alpha_1'$ and $\alpha_2'$ are MI. 
A similar argument can be made with the expression \eqref{eqn:ev2}.  
\end{remark}
\par
With the help of Lemma \ref{lem:DE}, FPF's \eqref{eqn:saito} and \eqref{eqn:TT} 
in Propositions \ref{prop:saito} and \ref{prop:TT} often make it possible to 
determine the rational functions $P(w)$ and $Q(w)$ explicitly. 
One more piece toward this calculation is to know how the map $f : X \to X$ 
permutes the $(-2)$-curves in $X$. 
As is remarked at the end of \S\ref{sec:hgm}, this can be done by calculating 
the action of $\tilde{A}$ on the simple system $\vD_{\rb}$ explicitly. 
Without doing so, however, it is sometimes feasible to get this information by 
looking at $\tilde{\varphi}_1(z)$ only.  
Recall from \eqref{eqn:pic-char} that $\tilde{\varphi}_1(z)$ is the characteristic 
polynomial of $\tilde{A}|\Pic$, so $\tilde{\varphi}_1(z)$ must be divisible by the 
characteristic polynomial $\chi(z)$ of $\tilde{A}|\Span \, \vD_{\rb}$. 
Thus the shape of $\tilde{\varphi}_1(z)$ constrains that of $\chi(z)$ and hence 
the way in which $\tilde{A}$ acts on $\vD_{\rb}$ to some extent or 
fully in some cases. 
Putting all these ingredients together, we are able to establish, for example,  
the following result.   
\begin{table}[h]
\centerline{
\begin{tabular}{ccccccccc}
\hline
     &            &           &                 &               &       &          &             &                    \\[-3mm]   
$\#$ & $\rho$ & $S(z)$ & $C(z)$ & $\psi(z)$ & ST & Dynkin & $\tilde{\varphi}_1(z)$ & $\Tr \, \tilde{A}$ \\
\hline
$1$ & $2$ & $\rS_1^{(20)}$ & $1$ & $\rS_1^{(10)} \rC_{21}$ & $\tau_7$ & $\rA_1$ & $\rC_1\rC_2$ & $1$ \\
$2$ & $4$ & $\rS_{22}^{(18)}$ & $\rC_4$ & $\rS_1^{(6)} \rC_{48}$ & $\tau_4$ & $\rA_1^{\oplus 2}$ & $\rC_1\rC_2 \rC_4$ & $-1$ \\
$3$ & $6$ & $\rS_5^{(16)}$ & $\rC_3\rC_4$ & $\rS_2^{(10)} \rC_{42}$ & $\tau_6$ & $\emptyset$ & $\rC_1\rC_2\rC_3\rC_4$ & $-1$ \\
$4$ & $8$ & $\rS_1^{(14)}$ & $\rC_{14}$ & $\rS_4^{(14)} \rC_{24}$ & $\tau_1$ & $\rE_8$ & $\rC_1^8$ & $8$ \\ 
$5$ & $10$ & $\rS_1^{(12)}$ & $\rC_{16}$ & $\rS_3^{(6)} \rC_{60}$ & $\tau_5$ & $\rD_9$ & $\rC_1^8 \rC_2^2$ & $7$ \\
$6$ & $12$ & $\rS_1^{(10)}$ & $\rC_4\rC_{16}$ & $\rS_1^{(6)}\rC_{40}$ & $\tau_2$ & $\rE_6 \oplus \rE_6$ & $\rC_1^4\rC_2^4\rC_4^2$ & $-1$ \\
$7$ & $14$ & $\rS_{16}^{(8)}$ & $\rC_3\rC_{12}\rC_{18}$ & $\rS_{12}^{(10)} \rC_{42}$ & $\tau_2$ & $\emptyset$ & 
$\rC_1\rC_2\rC_3\rC_{12}\rC_{18}$ & $-1$\\
$8$ & $16$ & $\rS_1^{(6)}$ & $\rC_4\rC_{26}$ & $\rS_{43}^{(22)}$ & $\tau_1$ & $\rD_{16}$ & $\rC_1^{16}$ & $16$ \\
$9$ & $18$ & $\rS_1^{(4)}$ & $\rC_8\rC_{12}\rC_{30}$ & see \eqref{eqn:psi18} & $\tau_1$ & $\rA_2^{\oplus 2} \oplus \rE_6 \oplus \rE_8$ & 
$\rC_1^{13}\rC_2^3 \rC_4$ & $11$ \\
\hline
\end{tabular}}
\caption{Some pairs $(\varphi, \psi)$ leading to K3 surface automorphisms with Siegel disks.}  
\label{tab:SD}
\end{table}
\begin{theorem} \label{thm:SD}
The pairs $(\varphi, \psi)$ in Table $\ref{tab:SD}$, which are obtained from 
Setups $\ref{setup1}$ and $\ref{setup2}$, lead to K3 surface automorphisms with 
Siegel disks, where $\rho$ is the Picard number and $\varphi(z) = S(z) \cdot C(z)$.    
\end{theorem}
{\it Proof}. 
Let $f : X \to X$ be the K3 surface automorphism lifted from $\tilde{A}$ and 
$\chi(z)$ be the characteristic polynomial of 
$f^*|\Span \, \vD_{\rb}(X) = \tilde{A}|\Span\, \vD_{\rb}$.  
Leaving entry $\#1$ in \S \ref{sec:P2} (see Theorem \ref{thm:P2}) 
we deal with the remaining entries. 
\par
For entry $\#2$ the map $f$ exchanges the two $\rA_1$-components of 
the exceptional set $\cE = \cE(X)$.  
For, otherwise, $\tilde{A}$ fixes the two simple roots in $\vD_{\rb}$, having 
at least two eigenvalues $1$, so $\chi(z)$ and hence $\tilde{\varphi}_1(z) = 
\rC_1(z) \cdot \rC_2(z) \cdot \rC_4(z)$ is divisible by $\rC_1(z)^2$, 
a contradiction. 
For entry $\#6$ a similar reasoning with $\tilde{\varphi}_1(z) = \rC_1(z)^4 \cdot 
\rC_2(z)^4 \cdot \rC_4(z)^2$ implies that $f$ exchanges the two 
$\rE_6$-components of $\cE$. 
Thus for these entries $f$ has no fixed points on $\cE$ and all fixed points  
of $f$ are isolated, that is, $N_f = 0$ and $\Fixi(f) = \Fix(f)$ in \eqref{eqn:saito}.     
This is also the case with entries $\#3$ and $\#7$ for which $\cE$ is empty. 
We have $\Tr f^*|H^2(X, \bC) = \Tr \tilde{A} = -1$ for these four entries. 
FPF \eqref{eqn:saito} then implies that $f$ admits a unique transverse 
fixed point $p \in X$. 
If the eigenvalues of $(d f)_p$ are expressed as in \eqref{eqn:ev1}, 
then FPF \eqref{eqn:TT} yields $1+\delta^{-1} = 
\{1-\delta^{\frac{1}{2}}(\alpha+\alpha^{-1}) + \delta\}^{-1}$ and hence    
$(\alpha+\alpha^{-1})^2 = P(\tau)$, where 
$$
P(w) := \dfrac{(w+1)^2}{w+2}. 
$$
For entry $\#2$ we have $0 < P(\tau_4) < 4$ and $P(\tau_8) > 4$; 
for entry $\#3$ we have $0 < P(\tau_6) < 4$ and $P(\tau_7) > 4$; 
for entry $\#6$ we have $0 < P(\tau_2) < 4$ and $P(\tau_4) > 4$;
for entry $\#7$ we have $0 < P(\tau_2) < 4$ and $P(\tau_3) > 4$. 
Therefore in these cases $p$ is the center of a Siegel disk by 
Lemma \ref{lem:mi1},(1)-(i).   
\par
For entry $\#4$ the exceptional set $\cE$ itself is the only connected 
component, which is of type $\rE_8$. 
We have $N_f = 1$ and $\mu(f, \cE) = 7$ from Lemma \ref{lem:DE},(1). 
So FPF \eqref{eqn:saito} with $\Tr f^*|H^2(X, \bC) = \Tr \tilde{A} = 8$ 
shows that $f$ has a unique transverse fixed point $p \in X \setminus \cE$. 
If the eigenvalues of $(d f)_p$ are expressed as \eqref{eqn:ev1}, then 
FPF \eqref{eqn:TT} reads 
$$
1+\delta^{-1} = \dfrac{1}{1-\delta^{\frac{1}{2}}(\alpha + \alpha^{-1}) + \delta} 
- \dfrac{\delta}{(1-\delta)^2}  \left( \dfrac{1}{1+\delta} 
+ \dfrac{1+\delta}{1+\delta+\delta^2} 
+ \dfrac{1+\delta + \delta^2 + \delta^3}{1+\delta + \delta^2 + \delta^3 +\delta^4}\right)  
+ \dfrac{1+\delta}{(1-\delta)^2}
$$
where the middle term in the RHS comes from \eqref{eqn:nupE} with $n = 8$.    
This equation gives $(\alpha+\alpha^{-1})^2 = P(\tau)$ with 
$$
P(w) := \dfrac{(w+2)(w^5 - 5 w^3 - w^2 + 5 w + 1)^2}{(w^5 + w^4 - 5 w^3 - 5 w^2 + 4 w + 3)^2}. 
$$
We observe $0 < P(\tau_1) < 4$ and $P(\tau_4) > 4$. 
Hence $p$ is the center of a Siegel disk by Lemma \ref{lem:mi1},(1)-(i).  
\par
For entry $\#5$ the exceptional set $\cE$ itself is the only connected 
component, which is of type $\rD_9$.  
The map $f$ acts on the dual graph $\varGamma$ of $\cE$ non-trivially. 
For, otherwise, $\tilde{A}$ fixes all simple roots in $\vD_{\rb}$, having at 
least nine eigenvalues $1$, so $\chi(z)$ and hence $\tilde{\varphi}_1(z) = 
\rC_1(z)^8 \cdot \rC_2(z)^2$ are divisible by $\rC_1(z)^9$, a contradiction. 
We have $N_f = 0$ and $\mu(f, \cE) = 8$ from Lemma \ref{lem:DE},(2). 
So FPF \eqref{eqn:saito} with $\Tr f^*|H^2(X, \bC) = \Tr \tilde{A} = 7$ 
shows that $f$ has a unique transverse fixed point $p \in X \setminus \cE$. 
If the eigenvalues of $(d f)_p$ are expressed as \eqref{eqn:ev1}, then 
FPF \eqref{eqn:TT} reads 
$$
1+\delta^{-1} = \dfrac{1}{1-\delta^{\frac{1}{2}}(\alpha + \alpha^{-1}) + \delta} 
+ \dfrac{1}{2(1+\delta)}  
+ \dfrac{1+\delta + \delta^2 + \delta^3 + \delta^4 + \delta^5 + \delta^6}{2(1+\delta^7)}, 
$$
where the last two terms in the RHS stem from \eqref{eqn:numD} with $n=9$.    
This equation gives $(\alpha+\alpha^{-1})^2 = P(\tau)$ with 
$$
P(w) := \dfrac{(w + 2) (w^4 - w^3 - 3 w^2 + w + 1)^2}{(w - 2)^2 (w + 1)^2 (w^2 + w - 1)^2}. 
$$
We observe $0 < P(\tau_5) < 4$ and $P(\tau_2) > 4$. 
Hence $p$ is the center of a Siegel disk by Lemma \ref{lem:mi1},(1)-(i).  
\par
For entry $\#8$ the exceptional set $\cE$ itself is the only connected 
component, which is of type $\rD_{16}$.  
The map $f$ acts on the dual graph $\varGamma$ of $\cE$ trivially, 
because $\tilde{\varphi}_1(z) = \chi(z) =\rC_1(z)^{16}$.  
We have $N_f = 1$ and $\mu(f, \cE) = 15$ from Lemma \ref{lem:DE},(1). 
So FPF \eqref{eqn:saito} with $\Tr f^*|H^2(X, \bC) = \Tr \tilde{A} = 16$ 
shows that $f$ has a unique transverse fixed point $p \in X \setminus \cE$. 
If the eigenvalues of $(d f)_p$ are expressed as \eqref{eqn:ev1}, then 
FPF \eqref{eqn:TT} reads 
$$
1+\delta^{-1} = \dfrac{1}{1-\delta^{\frac{1}{2}}(\alpha + \alpha^{-1}) + \delta} 
-\dfrac{\delta}{(1-\delta)^2} \left( \dfrac{2}{1+\delta} 
+ \dfrac{1+ \delta + \cdots + \delta^{12}}{1+\delta+ \cdots +\delta^{13}}\right)  
+\dfrac{1+\delta}{(1-\delta)^2},  
$$
where the middle term in the RHS stems from \eqref{eqn:nupD} with $n=16$.    
This equation gives $(\alpha+\alpha^{-1})^2 = P(\tau)$ with 
$$
P(w) := \dfrac{(w + 2) 
(w^8 - 2 w^7 - 6 w^6 + 11 w^5 + 11 w^4 - 16 w^3 - 7 w^2 + 5 w + 1)^2}{(w^3 - 3 w - 1)^2 
(w^5 - w^4 - 5 w^3 + 4 w^2 + 5 w - 3)^2}. 
$$
We observe $0 < P(\tau_1) < 4$ and $P(\tau_2) > 4$. 
Hence $p$ is the center of a Siegel disk by Lemma \ref{lem:mi1},(1)-(i).  
\par
For entry $\#9$ there are obvious decompositions $\cE = \cE(\rA_2^{\oplus 2}) 
\amalg \cE(\rE_6) \amalg \cE(\rE_8)$, 
$\vD_{\rb} = \vD_{\rb}(\rA_2^{\oplus}) \amalg \vD_{\rb}(\rE_6) \amalg \vD_{\rb}(\rE_8)$  
and $\Span \, \vD_{\rb} = \Span \, \vD_{\rb}(\rA_2^{\oplus}) \oplus \Span \, \vD_{\rb}(\rE_6) 
\oplus \Span \, \vD_{\rb}(\rE_8)$, preserved by $f = \tilde{A}$.   
Let $\chi(z) = \chi_1(z) \cdot \chi_2(z) \cdot \chi_3(z)$ be the corresponding decomposition 
of characteristic polynomials. 
As $\deg \chi(z) = 2 \times 2 + 6 + 8 = 13 \times 1 + 3 \times 1 + 2 = 
\deg \tilde{\varphi}_1(z)$, we have $\chi(z) = \tilde{\varphi}_1(z) = 
\rC_1(z)^{13} \cdot \rC_2(z)^3 \cdot \rC_4(z)$. 
Since $\tilde{A}$ acts on $\vD_{\rb}(\rE_8)$ trivially, we have $\chi_3(z) = \rC_1(z)^8$ and 
$\chi_1(z) \cdot \chi_2(z) = \rC_1(z)^5 \cdot \rC_2(z)^3 \cdot \rC_4(z)$. 
If $\tilde{A}$ acts on $\vD_{\rb}(\rE_6)$ trivially then $\chi_2(z) = \rC_1(z)^6$, which is absurd. 
Thus $\tilde{A}$ acts on $\vD_{\rb}(\rE_6)$ non-trivially, so that $\chi_2(z) = \rC_1(z)^4 \cdot 
\rC_2(z)^2$ and $\chi_1(z) = \rC_1(z) \cdot \rC_2(z) \cdot \rC_4(z)$. 
Consider further decompositions $\cE(\rA_2^{\oplus 2}) = \cE^+(\rA_2) \amalg 
\cE^-(\rA_2)$ and $\vD_{\rb}(\rA_2^{\oplus 2}) = \vD_{\rb}^+(\rA_2) \amalg \vD_{\rb}^-(\rA_2)$. 
Note that $\vD_{\rb}^{\pm}(\rA_2)$ are either preserved or permuted by $\tilde{A}$. 
In the former case we have a decomposition $\chi_1(z) = \chi_1^+(z) \cdot \chi_1^-(z)$ 
with each factor being either $\rC_1(z)^2$ or $\rC_1(z) \cdot \rC_2(z)$. 
This is impossible, so the latter is actually the case. 
Therefore, $f$ permutes $\cE^{\pm}(\rA_2)$, acts 
non-trivially on $\cE(\rE_6)$ and trivially on $\cE(\rE_8)$. 
\par
We then have $N_f = 1$ and $\mu(f, \cE) = \mu(f, \cE(\rE_6)) + \mu(f, \cE(\rE_8)) = 3 + 7 = 10$ 
from Lemma \ref{lem:DE},(1)-(2).   
So FPF \eqref{eqn:saito} with $\Tr f^*|H^2(X, \bC) = \Tr \tilde{A} = 11$ 
shows that $f$ has a unique transverse fixed point $p \in X \setminus \cE$. 
If the eigenvalues of $(d f)_p$ are expressed as \eqref{eqn:ev1}, then 
FPF \eqref{eqn:TT} can be represented as  
\begin{equation*}
\begin{split}
1+\delta^{-1} 
&= \dfrac{1}{1-\delta^{\frac{1}{2}} (\alpha + \alpha^{-1}) + \delta } 
+ \left\{ \dfrac{1}{2(1+\delta)} + \dfrac{1+\delta}{2(1+\delta^2)} \right\} 
\\[1mm]
& \phantom{=} - \dfrac{\delta}{(1-\delta)^2}  
\left( 
 \dfrac{1}{1+\delta} + \dfrac{1+\delta}{1+\delta+\delta^2} 
+ \dfrac{1+\delta + \delta^2 + \delta^3}{1+\delta + \delta^2 + \delta^3 +\delta^4} \right)  
+ \dfrac{1+\delta}{(1-\delta)^2},  
\end{split}
\end{equation*}
where the second and third terms in the RHS come from \eqref{eqn:numE} and 
\eqref{eqn:nupE} with $n=8$ respectively.     
This equation leads to $(\alpha+\alpha^{-1})^2 = P(\tau)$ with the rational function 
$$
P(w) := 
\dfrac{(w + 2) (w^3 - 4 w -  2)^2 (w^3 - w^2 - 2 w + 1)^2}{(w^2 - 2)^2 (w^4 - 4 w^2 - w + 1)^2}. 
$$
For $\tau = \tau_1$ we observe $0 < P(\tau) < 4$. 
Using the fact that $\tau$ has minimal polynomial $\ST^{(4)}_1(w) = w^2-w-3$,  
we can show that $P(\tau)$ has minimal polynomial $27 w^2- 11 w + 1$, 
which is not monic, so that $P(\tau)$ is not an algebraic integer.  
Therefore $p$ is the center of a Siegel disk by Lemma \ref{lem:mi1},(1)-(ii).                        
\hfill $\Box$
\begin{remark} \label{rem:SD}  
Table \ref{tab:SD} is just for the sake of illustration, providing only one example 
for each Picard number $\rho = 2, 4, 6, \dots, 18$. 
In fact there are much more pairs $(\varphi, \psi)$ leading to Siegel disks. 
As for examples with $\rho = 0$ we refer to McMullen \cite[Table 4]{McMullen1}.   
More examples in this case can be found in \cite[Tables 8.2, 8.3, 8.4]{IT}, 
which are constructed by the method of hypergeometric groups with $\psi(z)$ being 
an unramified Salem polynomial of degree $22$ and the matrix $B$, in place of $A$, 
playing the role of a Hodge isometry. 
\end{remark}
\section{Picard Number 2} \label{sec:P2}
Let $\lambda^{(20)}_1 \approx 1.2326135$ be the smallest Salem number of degree 
$20$, whose minimal polynomial is given by     
$$
\rS^{(20)}_1(z) = z^{20}-z^{19}-z^{15}+z^{14}-z^{11}+z^{10}-z^9+z^6-z^5-z+1. 
$$
A computer enumeration shows that the solutions to Setup \ref{setup1} with 
$S(z) = \rS^{(20)}_1(z)$ are given as in Table \ref{tab:S20-1}, where the meaning 
of the last S/H column becomes clear after Theorem \ref{thm:P2} is stated.    
\begin{table}[h]
\centerline{
\begin{tabular}{cccccclc|c}
\hline
          &           &           &           &           &           &                                &                       &       \\[-3mm]
$S(z)$ & $C(z)$ & $s(z)$ & $c(z)$ & $\ST$ & Dynkin & $\tilde{\varphi}_1(z)$ & $\Tr \tilde{A}$ & S/H \\ 
\hline
$\rS^{(20)}_{1}$ & $1$ & $\rS^{(10)}_{1}$ & $\rC_{21}$ & $\tau_7$ & $\rA_1$ & $\rC_1 \rC_2$ & $1$ & HS \\
$\rS^{(20)}_{1}$ & $1$ & $\rS^{(10)}_{2}$ & $\rC_{12} \rC_{20}$ & $\tau_6$ & $\rA_1$ & $\rC_1 \rC_2$ & $1$ & SS \\
$\rS^{(20)}_{1}$ & $1$ & $\rS^{(10)}_{9}$ & $\rC_{12} \rC_{24}$ & $\tau_9$ & $\rA_1$ & $\rC_1 \rC_2$ & $1$ & SH \\
$\rS^{(20)}_{1}$ & $1$ & $\rS^{(10)}_{9}$ & $\rC_{12} \rC_{30}$ & $\tau_6$ & $\rA_1$ & $\rC_1 \rC_2$ & $1$ & SS \\
$\rS^{(20)}_{1}$ & $1$ & $\rS^{(14)}_{1}$ & $\rC_{20}$ & $\tau_5$ & $\rA_1$ & $\rC_1 \rC_2$ & $1$ & SS \\
$\rS^{(20)}_{1}$ & $1$ & $\rS^{(14)}_{12}$ & $\rC_{20}$ & $\tau_3$ & $\rA_1$ & $\rC_1 \rC_2$ & $1$ & HS \\
$\rS^{(20)}_{1}$ & $1$ & $\rS^{(14)}_{12}$ & $\rC_{30}$ & $\tau_4$ & $\rA_1$ & $\rC_1 \rC_2$ & $1$ & SS \\
$\rS^{(20)}_{1}$ & $1$ & $\rS^{(14)}_{26}$ & $\rC_{30}$ & $\tau_7$ & $\rA_1$ & $\rC_1 \rC_2$ & $1$ & HS \\
$\rS^{(20)}_{1}$ & $1$ & $\rS^{(18)}_{4}$ & $\rC_{12}$ & $\tau_1$ & $\rA_1$ & $\rC_1 \rC_2$ & $1$ & SS \\
$\rS^{(20)}_{1}$ & $1$ & $\rS^{(18)}_{7}$ & $\rC_{12}$ & $\tau_6$ & $\rA_1$ & $\rC_1 \rC_2$ & $1$ & SS \\
$\rS^{(20)}_{1}$ & $1$ & $\rS^{(18)}_{16}$ & $\rC_{12}$ & $\tau_3$ & $\rA_1$ & $\rC_1 \rC_2$ & $1$ & HS \\
$\rS^{(20)}_{1}$ & 1$$ & $\rS^{(18)}_{32}$ & $\rC_{12}$ & $\tau_3$ & $\rA_1$ & $\rC_1 \rC_2$ & $1$ & HS \\
$\rS^{(20)}_{1}$ & $1$ & $\rS^{(22)}_{1}$ & $1$ & $\tau_5$ & $\rA_1$ & $\rC_1 \rC_2$ & $1$ & SS \\
$\rS^{(20)}_{1}$ & $1$ & $\rS^{(22)}_{3}$ & $1$ & $\tau_1$ & $\rA_1$ & $\rC_1 \rC_2$ & $1$ & SS \\
$\rS^{(20)}_{1}$ & $1$ & $\rS^{(22)}_{5}$ & $1$ & $\tau_3$ & $\rA_1$ & $\rC_1 \rC_2$ & $1$ & HS \\ 
\hline
\end{tabular}}
\caption{Picard number $\rho = 20$ (Setup \ref{setup1}).} 
\label{tab:S20-1}
\end{table}
Table \ref{tab:S20-1} then leads us to consider any K3 surface automorphism  
$f : X \to X$ such that    
\begin{itemize}
\setlength{\itemsep}{-1pt}
\item $X$ has Picard number $\rho(X) = 2$ and exceptional set $\cE(X)$ 
of Dynkin type $\rA_1$, 
\item $f$ has entropy $h(f) = \log \lambda^{(20)}_1$ and special eigenvalue 
$\delta = \delta(f)$ conjugate to $\lambda^{(20)}_1$, 
\item $f^*|\Pic(X)$ has characteristic polynomial 
$\tilde{\varphi}_1(z) = \rC_1(z) \cdot \rC_2(z) = (z-1)(z+1)$.  
\end{itemize} 
We remark that $\cE(X)$ consists of only one $(-2)$-curve $E \cong \bP^1$ and 
the special trace $\tau := \delta + \delta^{-1}$ is among the roots $\tau_1, \dots, \tau_9$ 
of the trace polynomial $\ST_1^{(20)}(w)$ such that $2 > \tau_1 > \dots > \tau_9 > -2$.  
\begin{theorem} \label{thm:P2} 
The map $f$ has exactly three fixed points in $X$ consisting of a pair $p_{\pm} \in E$ 
and a single point $p \in X \setminus E$. 
Each of them is either the center of a Siegel disks $(\mathrm{S})$ or a hyperbolic 
fixed point $(\mathrm{H})$, with $p_{\pm}$ being in the same case.        
How this dichotomy occurs is shown in Table $\ref{tab:P2}$ for each value of the 
special trace $\tau$.       
\begin{table}[hh] 
\centerline{
\begin{tabular}{l|ccccccccc}
\hline
ST & $\tau_1$ & $\tau_2$ & $\tau_3$ & $\tau_4$ & $\tau_5$ & $\tau_6$ & 
$\tau_7$ & $\tau_8$ & $\tau_9$ \\[1mm]
\hline
 & & & & & & & & &  \\[-3mm]
$p_{\pm}$ & S & S & H & S & S & S  & H & H & S \\[1mm]
$p$         & S & S & S & S & S  & S & S  & S & H \\[1mm]
\hline
\end{tabular}}
\caption{Center of a Siegel disk (S) or a hyperbolic fixed point (H).} 
\label{tab:P2} 
\end{table}
\end{theorem}
\par
Applying Theorem \ref{thm:P2} to the entries of Table \ref{tab:S20-1}, we obtain 
the S/H column in it, where for example HS and SS mean that $(p_{\pm}, p)$ is of 
types $(\mathrm{H}, \mathrm{S})$ and $(\mathrm{S}, \mathrm{S})$ respectively. 
Notice that all $\tau_1, \dots, \tau_9$ but $\tau_2$ and $\tau_8$ appear as 
special traces. 
Entry $\#1$ of Table \ref{tab:SD} is just the first entry of Table \ref{tab:S20-1}, 
thus the proof of Theorem \ref{thm:SD} is completed when Theorem \ref{thm:P2} 
is established.  
The rest of this section is devoted to the proof of Theorem \ref{thm:P2}. 
\par
In general if $F$ is a linear endomorphism with characteristic polynomial 
$\varphi(z)$, then by the relation between the generating function for power 
sums and that for elementary symmetric polynomials we have   
$$
\Tr(F^n) = \mbox{the coefficient of $z^n$ in the 
Maclaurin expansion of $- z \dfrac{d}{d z} \log \varphi^{\dagger}(z)$}  
\quad \mbox{for any} \quad n \ge 1,   
$$
where $\varphi^{\dagger}(z)$ is the reciprocal to $\varphi(z)$. 
Currently, $F = f^*|H^2(X, \bZ)$ is the induced map on middle cohomology 
group, having anti-palindromic characteristic polynomial 
$\varphi(z) = (z-1)(z+1) \, \rS^{(20)}_1(z)$. 
The above formula tells us that $\Tr(F^n) = 
\underline{1}, 3, \underline{1}, 3, 6, 3, \underline{1}, 3, \dots$ for 
$n = \underline{1}, 2, \underline{3}, 4, 5, 6, \underline{7}, 8, \dots$ respectively. 
In particular we notice   
\begin{equation} \label{eqn:Tr1}
\Tr(F) = \Tr(F^3) = \Tr(F^7) = 1.  
\end{equation}
This observation leads us to consider the map $f$ together with its third 
and seventh iterates $f^3$ and $f^7$. 
\par
For the M\"{o}bius transformation $f_E := f|_E$ there are four possibilities:  
\begin{enumerate}[(i)]
\setlength{\itemsep}{-1pt}
\item $f_E$ has two distinct fixed points $p_{\pm} \in E$ such that  
$(d f_E)_{p_{\pm}} = \beta^{\pm 1} \in \bC^{\times}$ with $\beta \neq 1$, 
$\delta^{\pm 1}$,    
\item $f_E$ has a unique fixed point $p_0 \in E$, in which case $(d f_E)_{p_0} = 1$ and 
$p_0 \in \Fixe_{\rI}(f)$,   
\item $f_E$ is an identity transformation, that is, $E$ is a fixed curve of $f$,  
\item $f_E$ has two distinct fixed points $p_{\pm} \in E$ such that 
$(d f_E)_{p_{\pm}} = \delta^{\pm 1}$, in which case $p_+ \in \Fixe_{\rII}(f)$.    
\end{enumerate} 
In case (i) the eigenvalues of $(d f)_{p_{\pm}}$ are $\beta^{\pm1}$ and 
$\delta \beta^{\mp1}$ as in \eqref{eqn:ev2}, so $p_{\pm}$ are transverse 
fixed points of $f$.  
\begin{lemma} \label{lem:case-iv} 
Case $(\mathrm{iv})$ does not occur. 
\end{lemma}
{\it Proof}. 
FPF \eqref{eqn:saito} is combined with equations \eqref{eqn:Tr1} and $N_f = 0$ 
to yield 
$$
3 = 2 + 1 \le \mu_{p_+}(f) + \mu_{p_-}(f) + \sum_{p \neq p_{\pm}} \mu_p(f) = 3, 
$$
where the sum is taken over all $p \in \Fixi(f)$ such that $p \neq p_{\pm}$. 
This shows that $\mu_{p_+}(f) = 2$, $\mu_{p_-}(f) = 1$ and $f$ 
has no other fixed points.  
The same is true for $f^3$. 
Since $(d f)_{p_-}$ has eigenvalues $\delta^{-1}$ and $\delta^2$, FPF 
\eqref{eqn:TT} and formula \eqref{eqn:iterate} in Theorem 
\ref{thm:iterate} for $n = 1$, $3$ lead to a system of equations,   
\begin{align*}
1+ \delta^{-1} 
&= \dfrac{1}{(1-\delta^{-1})(1-\delta^2)} + \dfrac{2 \delta + \theta}{(1-\delta)^2}, 
\\[1mm]
1+ \delta^{-3} 
&= \dfrac{1}{(1-\delta^{-3})(1-\delta^6)} + 
\dfrac{2 +4 \delta^3 + (1+\delta+\delta^2) \theta}{3 (1-\delta^3)^2}. 
\end{align*}
Eliminating $\theta$ from it we obtain an algebraic equation $(1+\delta) 
(3 + 5 \delta^2 - 2 \delta^3 + 9 \delta^4 - 2 \delta^5 + 5 \delta^6 + 3 \delta^8) = 0$ 
for $\delta$.  
This contradicts the fact that the minimal polynomial of $\delta$ is $\rS^{(20)}_1(z)$. 
Thus case (iv) cannot occur.  \hfill $\Box$ \par\medskip
Put $\sigma := \delta^{\frac{1}{2}} + \delta^{-\frac{1}{2}}$ with branch 
$\Re(\delta^{\frac{1}{2}}) > 0$. 
Note that $\sigma = \sqrt{\tau + 2} > 0$. 
\begin{lemma} \label{lem:case-i} 
In cases $(\mathrm{i})$, $(\mathrm{ii})$ and $(\mathrm{iii})$ the map $f$ has 
a unique fixed point $p \in X \setminus E$, which is transverse. 
Let $\delta^{\frac{1}{2}} \alpha^{\pm1} \in \bC^{\times}$ 
be the eigenvalues of $(d f)_p$ as in \eqref{eqn:ev1} and put 
$A := \alpha + \alpha^{-1}$ and $B := \beta + \beta^{-1}$, where by 
convention we understand that $\beta := 1$ and $B := 2$ in cases $(\mathrm{ii})$ 
and $(\mathrm{iii})$. 
Then $A$ and $B$ satisfy the equation    
\begin{equation} \label{eqn:AB1}
\sigma = \dfrac{\sigma}{\tau-B} + \dfrac{1}{\sigma -A},  
\end{equation}
where $(\tau-A)(\sigma - B)$ does not vanish.   
In terms of $B$ the number $A$ is expressed as  
\begin{equation} \label{eqn:A} 
A = \dfrac{(\tau+1) B + 2 - \tau^2}{\sigma(B+1-\tau)} \qquad 
\mbox{with} \quad B+1-\tau \neq 0. 
\end{equation}
\end{lemma}
{\it Proof}. 
In case (i) we have $\Fixi(f) \cap E = \{ p_{\pm}\}$, $\mu_{p_{\pm}}(f) = 1$ and 
$N_f = 0$. 
In case (ii) we have $\Fixi(f) \cap E = \{ p_0 \}$,  $\mu_{p_0}(f) = 2$ and 
$N_f = 0$. 
In case (iii) we have $\Fixi(f) \cap E = \emptyset$ and $N_f =1$.  
In any case FPF \eqref{eqn:saito} together with \eqref{eqn:Tr1} 
shows that $f$ has a unique fixed point $p \in X \setminus E$, 
which is simple, i.e. transverse. 
\par
In case (i), since $p_{\pm}$ are transverse fixed point of $f$, 
formula \eqref{eqn:h-index} gives   
\begin{equation*}  \label{eqn:ppm}
\nu_{p_+}(f) + \nu_{p_-}(f) = \dfrac{1}{(1-\beta)(1-\delta \beta^{-1})} + 
\dfrac{1}{(1-\beta^{-1})(1-\delta \beta)} = \dfrac{1+\delta}{1 -\delta B + \delta^2},  
\end{equation*}
so that FPF \eqref{eqn:TT} can be expressed as  
\begin{equation} \label{eqn:TT-i}
1+\delta^{-1} = \nu_{p_+}(f) + \nu_{p_-}(f) + \nu_p(f) = 
\dfrac{1+\delta}{1 -\delta B + \delta^2} + \dfrac{1}{1 - \delta^{\frac{1}{2}} A + \delta}, 
\end{equation}
which is multiplied by $\delta^{\frac{1}{2}}$ to yield equation \eqref{eqn:AB1}. 
In case (ii) the terms $\nu_{p_+}(f) + \nu_{p_-}(f)$ in \eqref{eqn:TT-i} should be 
replaced by $\nu_{p_0}(f)$, which is equal to $(1 + \delta)/(1-\delta)^2$ by 
formula \eqref{eqn:mu} in Theorem \ref{thm:mu}.   
In case (iii) those terms are not present, but instead a new $(1 + \delta)/(1-\delta)^2$ 
comes in due to the transition $N_f = 0 \mapsto N_f = 1$.   
In either case we have \eqref{eqn:TT-i} and hence \eqref{eqn:AB1} with convention 
$B = 2$. 
Note that $\tau-B \neq 0$ follows from the transversality of $p_{\pm}$ in 
case (i) and from $B = 2$ in cases (ii) and (iii), while $\sigma -A \neq 0$ follows 
from the transversality of $p$.  
\par
Equation \eqref{eqn:AB1} yields $\sigma(B+1-\tau) A = (\tau+1) B + 2 - \tau^2$. 
If $B+1-\tau = 0$, that is, $B = \tau-1$ then $0 = (\tau+1) B + 2 - \tau^2 
= 1$, which is impossible. 
Hence $B+1-\tau \neq 0$ and $A$ is expressed as \eqref{eqn:A}.  
\hfill $\Box$
\begin{lemma} \label{lem:delta37}
In case $(\mathrm{i})$ we have $\beta^n \neq \delta^{\pm n}$ for $n = 3$, $7$.  
\end{lemma}
{\it Proof}. If $\beta^3 = \delta^{\pm 3}$ then $p_{\pm}$, $p \in \Fixi(f^3)$ 
with $\mu_{p_{\pm}}(f^3) \ge 2$, $\mu_{p_{\mp}}(f^3) = 1$, $\mu_p(f^3) \ge 1$  
and $N_{f^3} = 0$, so FPF \eqref{eqn:saito} together with 
\eqref{eqn:Tr1} leads to a contradiction   
$4 = 2 + 1 + 1 \le \mu_{p_{\pm}}(f^3) + \mu_{p_{\mp}}(f^3) + \mu_p(f^3) \le 3$.  
If $\beta^7 = \delta^{\pm 7}$ then the same argument with $f^3$ replaced by 
$f^7$ yields a similar contradiction. 
\hfill $\Box$ \par\medskip
For $n \ge 1$ let $h_n(w) \in \bZ[w]$ be the polynomial such that 
$z^n + z^{-n} = h_n(w)$ for $w = z + z^{-1}$. 
We have 
$$
h_3(w) = w(w^2-3), \qquad h_7(w) = w(w^6-7 w^4+14 w^2-7). 
$$
\begin{lemma} \label{lem:AB37} 
In cases $(\mathrm{i})$, $(\mathrm{ii})$ and $(\mathrm{iii})$ the numbers 
$A$ and $B$ in Lemma $\ref{lem:case-i}$ satisfy two more equations   
\begin{equation} \label{eqn:AB37} 
h_n(\sigma) = \dfrac{h_n(\sigma)}{h_n(\tau) - h_n(B)} 
+ \dfrac{1}{h_n(\sigma) - h_n(A)},   \qquad n = 3, \, 7, 
\end{equation}
where all fractions appearing in \eqref{eqn:AB37} have nonzero denominators.  
\end{lemma}
{\it Proof}. For $n = 3$, $7$ we have $\Tr(f^n) = \Tr(f) = 1$ by \eqref{eqn:Tr1} 
and
\begin{itemize}
\setlength{\itemsep}{-1pt}
\item case (i) for $f$ with $\beta^n \neq 1$ leads to case (i) for $f^n$  
by Lemma \ref{lem:delta37}, 
\item case (i) for $f$ with $\beta^n = 1$ leads to case (iii) for $f^n$, 
in which $h_n(B) = \beta^n + \beta^{-n} = 2$,     
\item case (ii) for $f$ leads to case (ii) for $f^n$, in which  
$h_n(B) = h_n(2) = 2$,  
\item case (iii) for $f$ leads to case (iii) for $f^n$, in which 
$h_n(B) = h_n(2) = 2$.    
\end{itemize}
Thus Lemma \ref{lem:case-i} and its proof apply to $f^n$ in place of $f$. 
Equations \eqref{eqn:AB37} are obtained from \eqref{eqn:AB1} by replacing 
$f$ with $f^n$. 
This amounts to altering $\delta \mapsto \delta^n$, $\alpha \mapsto \alpha^n$, 
$\beta \mapsto \beta^n$ and so $\xi \mapsto h_n(\xi)$ for $\xi = \tau$, 
$\sigma$, $A$, $B$ in \eqref{eqn:AB1}. 
Here in cases (ii) and (iii) the convention in Lemma \ref{lem:case-i} takes the 
form $h_n(B) = 2$ for $f^n$, which is fulfilled.  
\hfill $\Box$
\begin{lemma} \label{lem:case-ii} 
Cases $(\mathrm{ii})$ and $(\mathrm{iii})$ do not occur, 
hence case $(\mathrm{i})$ actually occurs.    
\end{lemma}
{\it Proof}.  
Recall that we have $B = 2$ in cases (ii) and (iii). 
Substituting \eqref{eqn:A} with $B=2$ into \eqref{eqn:AB37} for $n = 3$, 
we find that $\tau$ satisfies the septic equation 
$\tau^7 - 3 \tau^6 - 9 \tau^5 + 17 \tau^4 + 39 \tau^3 - 7 \tau^2 - 50 \tau - 20 = 0$.  
This contradicts the fact that the minimal polynomial of $\tau$ is $\ST^{(20)}_1(w)$.  
Thus these cases cannot occur altogether. \hfill $\Box$
\begin{lemma} \label{lem:quartic} 
Let $A$ and $B$ be the numbers in Lemma $\ref{lem:case-i}$. 
Then we have 
\begin{subequations} \label{eqn:roots}
\begin{align}
B &= Q(\tau) := -(\tau+1)(\tau-2)(\tau^3-3 \tau +1), \label{eqn:rootsB} \\[1mm]
A^2 &= P(\tau) := \dfrac{(\tau^6-6\tau^4-\tau^3+10 \tau^2+3 \tau-4)^2}{(\tau+2) 
(\tau^2-3)^2(\tau^3-\tau^2-2\tau+1)^2}. \label{eqn:rootsA}
\end{align}
\end{subequations}
\end{lemma}
{\it Proof}. 
Substituting \eqref{eqn:A} into \eqref{eqn:AB37} we obtain two algebraic 
equations for $B$, which turn out to factor into  
$$
\{ B- Q(\tau) \} \, R_3(\tau; B) = 0, \qquad \{ B - Q(\tau) \} \, R_7(\tau; B) = 0, 
$$
over the number field $K := \bQ(\tau)$, where $R_3(\tau; x) \in K[x]$ and 
$R_7(\tau; x) \in K[x]$ are polynomials of degrees $3$ and $11$ respectively. 
Moreover $R_3(\tau; x)$ and $R_7(\tau; x)$ have no roots in common (consider 
their resultant).   
These facts are verified by Mathematica, which is capable of polynomial calculations 
over an algebraic number field. 
Thus we obtain equation \eqref{eqn:rootsB}.  
Substituting it into \eqref{eqn:A} yields  
$$
A = \dfrac{\tau^6-6\tau^4-\tau^3+10 \tau^2+3 \tau-4}{\sigma 
(\tau^2-3)(\tau^3-\tau^2-2\tau+1)},  
$$
which is squared to give equation \eqref{eqn:rootsA}, where the relation 
$\sigma^2 = \tau + 2$ is also used.  \hfill $\Box$ \par\medskip
{\it Proof of Theorem $\ref{thm:P2}$}. 
We observe that $|Q(\tau_j)| <2$ for $j = 1, 2, 4, 5, 6, 9$ and $|Q(\tau_j)| > 2$ for 
$j = 3, 7, 8$. 
Thus Lemma \ref{lem:mi2} together with Remark \ref{rem:mi} implies the second row 
in Table \ref{tab:P2}.  
Similarly we observe that $0 < P(\tau_j) < 4$ for $j = 1, \dots, 8$ and $P(\tau_9) > 4$. 
Hence Lemma \ref{lem:mi1} yields the third row in Table \ref{tab:P2}.  \hfill $\Box$
\begin{remark} \label{rem:A2} 
There are examples of K3 surface automorphisms $f : X \to X$ such that 
$\rho(X) = 12$, the exceptional set $\cE(X)$ is of type $\rA_2$, and $f$ has three 
Siegel disks with centers on $\cE(X)$ (see \cite[Remark 9.7]{IT}).  
\end{remark}
\appendix
\section{Table of Salem Trace Polynomials} \label{sec:gloss}
Let $\lambda_i^{(d)}$ be the $i$-th smallest Salem number of degree $d$ and 
$\rS_i^{(d)}(z)$ be its minimal polynomial. 
Here is a list of all Salem polynomials $\rS_i^{(d)}(z)$ that appear explicitly in this 
article as the Salem factor $S(z)$ of the polynomial $\varphi(z)$. 
They are presented in terms of their trace polynomials $\ST_i^{(d)}(w)$.  
For each of them numerical computations and symbolic manipulations of the roots 
$\tau_0, \tau_1, \dots, \tau_{d/2-1}$ in \eqref{eqn:trace} can be carried out by 
using these data.  
\begin{equation*}
\begin{split}
\ST_1^{(4)}(w) &= w^2 - w - 3, \\
\ST_1^{(6)}(w) &= w^3 - 4 w - 1, \\
\ST_1^{(8)}(w) &= w^4 - 4 w^2 - w + 1, \\
\ST_2^{(8)}(w) &= w^4 - w^3 - 3 w^2 + w + 1, \\
\ST_{15}^{(8)}(w) &= w^4 - 2 w^3 - 4 w^2 + 7 w + 1, \\
\ST_{16}^{(8)}(w) &= w^4 - 5 w^2 - 2 w + 1, \\
\ST_1^{(10)}(w) &=  w^5 + w^4 - 5 w^3 - 5 w^2 + 4 w + 3, \\
\ST_1^{(12)}(w) &= w^6 - w^5 - 5 w^4 + 4 w^3 + 5 w^2 - 2 w - 1, \\
\ST_1^{(14)}(w) &= w^7 - 7 w^5 - w^4 + 13 w^3 + 4 w^2 - 4 w - 1, \\
\ST_1^{(16)}(w) &= w^8 - w^7 - 8 w^6 + 7 w^5 + 20 w^4 - 14 w^3 - 16 w^2 + 7 w + 1, \\
\ST_2^{(16)}(w) &= w^8 + w^7 - 8 w^6 - 8 w^5 + 19 w^4 + 18 w^3 - 13 w^2 - 10 w + 1, \\
\ST_3^{(16)}(w) &= w^8 - 8 w^6 - w^5 + 20 w^4 + 4 w^3 - 16 w^2 - 3 w + 2, \\
\ST_4^{(16)}(w) &= w^8 - w^7 - 8 w^6 + 7 w^5 + 20 w^4 - 14 w^3 - 17 w^2 + 7 w + 4, \\
\ST_5^{(16)}(w) &= w^8 - 9 w^6 - w^5 + 26 w^4 + 5 w^3 - 25 w^2 - 5 w + 4, \\
\ST_{22}^{(18)}(w) &= w^9 + w^8 - 10 w^7 - 11 w^6 + 32 w^5 + 38 w^4 - 33 w^3 - 42 w^2 + 4 w + 7, \\
\ST_1^{(20)}(w) &= w^{10} - w^9 - 10 w^8 + 9 w^7 + 35 w^6 - 28 w^5 - 49 w^4 + 35 w^3  + 21 w^2 -15 w + 1. 
\end{split}
\end{equation*}
\par\vspace{2mm} \noindent
{\bf Acknowledgments}. 
This work was supported by JSPS KAKENHI Grant Numbers JP19K03575, JP21J20107.    
 
\end{document}